\documentclass[final,onefignum,onetabnum]{siamart220329}

\usepackage{amsfonts}
\usepackage{graphicx}
\usepackage{epstopdf}
\usepackage{algorithmic}
\usepackage{amsmath}
\usepackage{bm}
\usepackage{dutchcal} 
\usepackage{dsfont}
\usepackage{stmaryrd}
\usepackage{verbatim}
\usepackage{enumitem}

\ifpdf
\DeclareGraphicsExtensions{.eps,.pdf,.png,.jpg}
\else
\DeclareGraphicsExtensions{.eps}
\fi

\usepackage{xcolor}
\definecolor{darkgreen}{rgb}{0.1,0.6,0.1}

\allowdisplaybreaks

\newsiamthm{remark}{Remark}
\newsiamthm{hypothesis}{Hypothesis}
\newsiamthm{schemedefinition}{Definition}
\crefname{hypothesis}{Hypothesis}{Hypotheses}
\crefname{schemedefinition}{Definition}{Definitions}
\newsiamthm{claim}{Claim}
\newsiamthm{schemeproperty}{Property}
\newsiamthm{schemetheorem}{Theorem}
\newsiamthm{schemeproposition}{Proposition}
\crefname{schemeproperty}{Property}{Properties}
\crefname{schemetheorem}{Theorem}{Theorems}
\crefname{schemeproposition}{Proposition}{Propositions}

\headers{Involution-Preserving ALE schemes}{S. Cardonna, E. Bernardelli and W. Boscheri}

\title{Involution-Preserving schemes through constrained mesh motion: a novel Arbitrary-Lagrangian-Eulerian framework}

\author{
	Sacha Cardonna\thanks{Institut Montpelliérain Alexander Grothendieck UMR 5149, 34090 Montpellier, France (\email{sacha.cardonna@umontpellier.fr})}
	\and
	Elena Bernardelli\thanks{Department of Computer Science, University of Verona, 37134 Verona, Italy (\email{elena.bernardelli@univr.it})}
	\and
	Walter Boscheri\thanks{Laboratoire de Mathématiques UMR 5127 CNRS, Universit{\'e} Savoie Mont Blanc, 73376 Le Bourget du Lac, France (\email{walter.boscheri@univ-smb.fr}) and Department of Mathematics and Computer Science, University of Ferrara, via Machiavelli 30, 44121 Ferrara, Italy}
}

\makeatletter
\AtBeginDocument{
	\def\refstepcounter@optarg[#1]#2{%
		\cref@old@refstepcounter{#2}%
		\cref@constructprefix{#2}{\cref@result}%
		\@ifundefined{cref@#1@alias}%
		{\def\@tempa{#1}}%
		{\def\@tempa{\csname cref@#1@alias\endcsname}}%
		\protected@edef\cref@currentlabel{%
			[\@tempa][\arabic{#2}][\cref@result]%
			\csname p@#2\endcsname\csname the#2\endcsname}}}      
\makeatother

\begin{document}

\maketitle
 
\begin{abstract} 
We introduce a novel structure-preserving finite volume method that exploits mesh motion to enforce differential involutions arising in hyperbolic systems of conservation laws. The key idea is to leverage the geometric degrees of freedom of an Arbitrary-Lagrangian–Eulerian (ALE) framework to construct a constrained mesh motion that is intrinsically compatible with the involution constraints embedded in the governing equations. The method is developed within a cell-centered formulation on moving quadrilateral meshes, where discrete divergence and curl operators are constructed from vertex-based differential operators. Particular care is devoted to the consistent evolution of the mesh geometry and the satisfaction of the Geometric Conservation Law (GCL).

The proposed framework is applied to two representative hyperbolic systems featuring distinct involution constraints: the ideal magnetohydrodynamics equations, endowed with the divergence-free condition on the magnetic field, and the hyperbolic Euler–Heat equations, subject to a curl-free constraint on the thermal impulse. Second-order accuracy in time is achieved through a Runge–Kutta discretization, with the geometric constraints consistently enforced at each stage. A comprehensive set of classical benchmark problems is presented to assess the accuracy and effectiveness of the method. The numerical results demonstrate that appropriately constrained mesh motion provides an effective mechanism for preserving differential involutions, establishing a novel connection between mesh geometry and the structure-preserving discretization of hyperbolic systems.

\end{abstract}

\begin{keywords}
structure-preserving schemes; involution preservation; continuum mechanics; moving unstructured meshes; constrained mesh motion
\end{keywords}

\begin{AMS}
  35L40, 
  65M08. 
\end{AMS}

\section{Introduction}\label{sec.intro}

Many systems arising in continuum mechanics are endowed with additional geometric and structural constraints, commonly referred to as involutions. Among these, divergence-free and curl-free conditions play a particularly important role. Such constraints arise in several physical models, 
and their preservation at the discrete level is particularly important, since violations may deteriorate the quality of the numerical solution, lead to numerical instabilities or finite-time blow-up \cite{chiocchetti2023exactly,dhaouadi2023structure}, and result in the loss of relevant physical properties \cite{balsara1999staggered,brackbill1985fluid,brackbill1980effect, bernardelli2026structure,Sidilkover2025Spurious}. Here, we consider two representative hyperbolic systems:
(i) the equations of magnetohydrodynamics (MHD), which describe the 
dynamics of electrically conducting fluids, such as plasma, through 
the coupling of fluid motion and electromagnetic fields; and
(ii) the first-order hyperbolic Euler--Heat model for heat-conducting compressible 
fluids introduced in \cite{DhaouadiGavrilyuk2024,BoscheriDhaouadi2026}, which extends the compressible Euler equations to incorporate heat conduction through 
an additional vector field, the thermal impulse $\bm{j}$. For the MHD equations, the relevant involution is the solenoidal constraint on the magnetic field, $\nabla\cdot\bm{B}=0$, which follows from Maxwell equations and expresses the absence of magnetic monopoles. For the hyperbolic Euler--Heat model, the evolution equation for the thermal impulse has a gradient structure, which implies the curl-free constraint $\nabla\times\bm{j}=\mathbf{0}$. If satisfied at the initial time, these conditions are preserved by the continuous evolution of the systems. Satisfying these structures at the discrete level of the numerical method remains an open challenge.

The numerical approximation of hyperbolic systems on moving computational meshes is naturally formulated within the Arbitrary Lagrangian--Eulerian (ALE) framework \cite{dumbser2013arbitrary, gaburro2026treatment,ShashkovKnupp2002,Loubere2010,Barlow2016}. In contrast to the Eulerian description, where the computational mesh remains fixed, and the Lagrangian description, where the mesh moves with the local fluid velocity, ALE methods allow the grid nodes to move with an arbitrary prescribed velocity. This additional degree of freedom combines the advantages of Eulerian and Lagrangian approaches while providing greater control over the evolution and quality of the computational mesh. In particular, by adapting the mesh motion to the underlying flow, ALE schemes can reduce the numerical dissipation associated with convective transport and improve the resolution of moving discontinuities, contact waves, and material interfaces. At the same time, the freedom in the choice of the mesh velocity can be exploited to limit excessive mesh deformation, which represents one of the main difficulties of purely Lagrangian approaches.
A fundamental requirement of numerical schemes on moving meshes is the satisfaction of the Geometric Conservation Law (GCL), which ensures that the evolution of the computational geometry is consistently coupled with the discrete governing equations. In particular, satisfaction of the GCL is closely related to the preservation of uniform flow under arbitrary mesh motion and is essential to avoid errors generated solely by the deformation of the computational grid. While the GCL guarantees consistency between the numerical solution and the evolving geometry, additional difficulties arise when the governing PDE system is subject to differential involutions. 

In this work, we introduce a novel involution-preserving ALE finite volume approach based on constrained mesh motion. The key idea is to exploit the freedom in the choice of the ALE mesh velocity and determine it through a constrained problem that accounts for the differential structure of the underlying physical system. The proposed method is formulated in a cell-centered finite volume framework on moving structured quadrilateral meshes. The discrete differential operators are constructed from corner-based geometrical quantities, providing a natural discrete representation of the divergence and curl operators on the evolving mesh. The computational mesh is consistently updated during the mesh evolution, while a second-order Runge--Kutta method is employed for the time discretization. The mesh motion is constrained with the aim of preserving the involution over time, according to the corresponding differential structure, either the divergence or the curl operator.

The outline of the paper as follows. In Section \ref{sec.model} we introduce the two systems of partial differential equations treated in this paper as representative examples. In Section \ref{sec.geom} we present the geometrical and discrete ingredients used by the novel approach, such as the corner vectors and the discrete operators to compute the involutions. Section \ref{sec.method} describes the fully discrete structure-preserving mesh motion. In Section \ref{sec.results}, a validation through classical benchmark test cases is shown. Conclusions and an outlook to further work are given in Section \ref{sec.concl}.

\section{Governing equations and continuous involutions}
\label{sec.model}

Let $\Omega\subset\mathbb{R}^2$ be a bounded domain with boundary $\partial \Omega$, let $\bm{x}=(x,y)\in\Omega$ denote the spatial coordinates, and let $t\in\mathbb{R}_+$ be the time coordinate. The governing equations can be written in the general form
\begin{eqnarray}
	\frac{\partial \bm{q}}{\partial t} + \nabla \cdot \mathds{F}(\bm{q}) = \bm{0},
	\label{eqn.pde}
\end{eqnarray}
where $\bm{q} \in \Omega_{\bm{q}}$ is the vector of conserved variables defined in the set of admissible states $\Omega_{\bm{q}}$, and $\mathds{F}(\bm{q})$ is the nonlinear flux tensor. Vectors are denoted by bold symbols and double-stroke symbols are used for tensors.



\subsection*{Ideal magnetohydrodynamics (MHD)}
\noindent The MHD equations can be cast into the formulation \eqref{eqn.pde} by defining the state vector and the flux tensor
\begin{align}
  \bm{q}
  =
  \begin{pmatrix}
    \rho\\ \rho\bm{u}\\ \rho E\\ \bm{B}
  \end{pmatrix}, \qquad \mathds{F}(\bm{q})
  =
  \begin{pmatrix}
  \rho\bm{u}\\[1mm]
  \rho\bm{u}\otimes\bm{u}
  +\left(p+\dfrac{\lVert\bm{B}\rVert^2}{8\pi}\right)\mathds{I}_2
  -\dfrac{1}{4\pi}\bm{B}\otimes\bm{B}\\[2mm]
  \left(\rho E+p+\dfrac{\lVert\bm{B}\rVert^2}{8\pi}\right)\bm{u}
  -\dfrac{1}{4\pi}\bm{B}(\bm{u}\cdot\bm{B})\\[2mm]
  \bm{B}\otimes\bm{u}-\bm{u}\otimes\bm{B}
  \end{pmatrix}.
  \label{eq:mhd_state}
\end{align}
Here, $\rho$ is the fluid density, $\bm{u}=(u,v)$ is the velocity, $p$ is the pressure, $\rho E$ is the total energy density, and $\bm{B}=(B_x,B_y)$ is the magnetic field. The system is closed with the ideal gas equation of state, thus the the total energy density is given by
\begin{align}
  \rho E
  =
  \rho e+\frac{1}{2}\rho\lVert\bm{u}\rVert^2
  +\frac{\lVert\bm{B}\rVert^2}{8\pi}, \qquad \rho e=\frac{p}{\gamma-1},
  \label{eq:mhd_energy}
\end{align}
with $\gamma>1$ representing the ratio of specific heats of the gas. The magnetic field is subject to the solenoidal constraint $\nabla\cdot\bm{B}=0$, which must be preserved over time. 

\subsection*{Hyperbolic Euler--Heat model (EH)}
\noindent The EH model is given by defining
\begin{align}
  \bm{q}
  =
  \begin{pmatrix}
    \rho\\ \rho\bm{u}\\ {\rho E}\\ \bm{j}
  \end{pmatrix}, \qquad \mathds{F}(\bm{q})
  =
  \begin{pmatrix}
  \rho\bm{u}\\[1mm]
  \rho\bm{u}\otimes\bm{u}
  +\left(p-\alpha(\rho)\lVert\bm{j}\rVert^2\right)\mathds{I}_2
  +\alpha(\rho)\bm{j}\otimes\bm{j}\\[1mm]
  \left({\rho E}+p-\alpha(\rho)\lVert\bm{j}\rVert^2\right)\bm{u}
  +\alpha(\rho)\phi\,\bm{j}\\[1mm]
  \phi\,\mathds{I}_2
  \end{pmatrix},
  \label{eq:eh_state}
\end{align}
with $\phi=\bm{j}\cdot\bm{u}+\theta$ and the thermal impulse $\bm{j}=(j_x,j_y)$.
The system is closed by the constitutive relations
\begin{align}
  \alpha(\rho)=\frac{\kappa^2}{\rho},
  \qquad
  {\rho E}
  =
  \frac{1}{2}\rho\lVert\bm{u}\rVert^2
  +\frac{p}{\gamma-1}
  +\frac{1}{2}\alpha(\rho)\lVert\bm{j}\rVert^2,
  \qquad
  \theta=\frac{p}{\rho c_v(\gamma-1)},
  \label{eq:eh_closure}
\end{align}
with $\kappa\geq0$, $c_v>0$, and $\gamma>1$. For this choice of $\alpha$, the general isotropic impulse contribution $\tfrac{1}{2}(\rho\alpha'(\rho)-\alpha(\rho)) \lVert\bm{j}\rVert^2$ reduces to $-\alpha(\rho)\lVert\bm{j}\rVert^2$.
The thermal impulse obeys the involution $\nabla\times\bm{j}=\mathbf{0}$.

\begin{remark}
The complete EH model also contains the rotational term $[\nabla\bm{j}-(\nabla\bm{j})^\top]\bm{u}$, which vanishes on the curl-free manifold, which is considered here. Furthermore, a source term of the form $-\bm{j}/\tau$ might be added to the EH system with a constant relaxation parameter $\tau$, which still implies a curl-free manifold.
\end{remark}

\section{Arbitrary-Eulerian-Lagrangian (ALE) finite volume scheme} \label{sec.geom}

At time $t$, the domain is covered by a conforming tessellation $\mathcal{T}_h(t)=\{\omega_c(t)\}_{c=1}^{N_c}$ made of $N_c$ non-overlapping time-dependent control volumes $\omega_c(t)$ with boundary $\partial\omega_c(t)$ and outward pointing normal vector $\bm{n}$. Let $N_p$ be the total number of vertices. We define the complete, boundary, and interior nodal index sets by
\begin{align}
  \mathcal{N}_h&:=\llbracket 1,N_p\rrbracket,
  &
  \mathcal{N}_h^\partial&:=\{p\in\mathcal{N}_h\mid\bm{x}_p(t)\in\partial\Omega\},
  &
  \mathcal{N}_h^\circ&:=\mathcal{N}_h\setminus\mathcal{N}_h^\partial.
  \label{eq:nodal_index_sets}
\end{align}
The set of vertices of $\omega_c$ is denoted by $\mathcal{P}(c)\subset\mathcal{N}_h$, while $\mathcal{C}(p)$ is the set of cells sharing vertex $p$; $\sharp\mathcal{P}(c)$ denotes the cardinality of $\mathcal{P}(c)$. The position and grid velocity of vertex $p\in\mathcal{N}_h$ are denoted by $\bm{x}_p(t)$ and $\bm{\mathcal{w}}_p(t)$, respectively. Let us introduce $\bm{X}:=(\bm{x}_p)_{p\in\mathcal{N}_h}$ to represent the complete vector of nodal coordinates. The mesh topology is kept fixed during the computation; only the vertex coordinates evolve.

\medskip
\noindent
\begin{schemedefinition}[Corner geometry]
\label{def:corner_geometry}
Let $p^-$ and $p^+$ be the previous and next vertices of $p$ along the counter-clockwise oriented boundary of $\omega_c$. The two edges adjacent to the corner $(p,c)$ have lengths $\ell_{pc}^-$ and $\ell_{pc}^+$ and outward unit normals $\bm{n}_{pc}^-$ and $\bm{n}_{pc}^+$. The associated corner vector is
\begin{align}
  \bm{\ell}_{pc}
  =
  \frac{1}{2} \left( \ell_{pc}^- \bm{n}_{pc}^- +
\ell_{pc}^+ \bm{n}_{pc}^+ \right).
  \label{eq:corner_vector}
\end{align}
It contains the two half-edge contributions attached to vertex $p$ and satisfies the discrete Gauss theorem
\begin{align}
  \sum_{p\in\mathcal{P}(c)}\bm{\ell}_{pc}=\bm{0}.
  \label{eq:corner_closure}
\end{align}
We further introduce the tangent corner vector obtained by a counter-clockwise rotation, that is
\begin{align}
  \bm{t}_{pc}
  :=
  \bm{\ell}_{pc}^{\perp}
  =
  \begin{pmatrix}
    -(\bm{\ell}_{pc})_2\\
    (\bm{\ell}_{pc})_1
  \end{pmatrix}.
  \label{eq:tangent_corner}
\end{align}
\end{schemedefinition}

\begin{schemedefinition}[Dual cell and subcells]
\label{def:nodal_patch}
The dual cell $\omega_p$ associated with a mesh vertex $p$ is constructed by connecting the barycenters of the cells sharing $p$ with the midpoints of the edges incident to $p$, thereby defining a vertex-centered
polygonal control volume surrounding $p$. Each primal cell is further subdivided into subcells $\omega_{pc}$,
each bounded by the cell barycenter, a vertex $p$, and the midpoints
of the two cell edges incident to $p$. The volume of the dual cell is then given by
\begin{align}
  |\omega_p|
  =
  \sum_{c\in\mathcal{C}(p)}|\omega_{pc}|,
  \label{eq:nodal_volume}
\end{align}
and the fundamental geometry identity \eqref{eq:corner_closure} holds true on the dual cell, that is
\begin{align}
	\sum_{c\in\mathcal{C}(p)}\bm{\ell}_{cp}=\bm{0},
	\label{eq:nodal_corner_closure}
\end{align}
with $\bm{\ell}_{cp}=-\bm{\ell}_{pc}$ being the corner normal vector of the dual cell.
\end{schemedefinition}

We design an ALE finite volume method closely related to the framework used in cell-centered Lagrangian schemes \cite{Maire2009,BoscheriGPRGCL}. Within each cell $\omega_c$, the solution is stored at time $t$ as cell-averages given by
\begin{align}
  \bm{q}_c
  =
  \frac{1}{|\omega_c|}
  \int_{\omega_c}\bm{q}(\bm{x})\,\mathrm{d}\bm{x}.
  \label{eq:cell_average}
\end{align}
The cell averages are transferred to the dual cells through the volume-weighted average
\begin{align}
  \bm{q}_p
  =
  \frac{1}{|\omega_p|}
  \sum_{c\in\mathcal{C}(p)}
  |\omega_{pc}|\,\bm{q}_c .
  \label{eq:nodal_state}
\end{align}
The semi-discrete ALE finite volume scheme writes
\begin{align}
\frac{\mathrm{d}}{\mathrm{d}t}
\left(|\omega_c|\bm{q}_c\right)
+
\sum_{p\in\mathcal{P}(c)}
\widehat{\mathds{F}}_p \, \bm{\ell}_{pc}
=\bm{0},
\label{eq:node_fv_scheme}
\end{align}
with the nodal ALE flux given by the following multidimensional Riemann solver 
\begin{align}
\widehat{\mathds{F}}_p
=
\mathds{F}_p
-\frac{1}{2}\sqrt{|\omega_p|}\,
a_p\,\mathds{G}_p, \qquad \mathds{F}_p = \mathds{F}(\bm{q}_p) -\bm{q}_p \otimes \bm{\mathcal{w}}_p. 
\label{eq:stabilized_nodal_flux}
\end{align}
The nodal gradient is computed as
\begin{align}
  \mathds{G}_p
  :=  -\frac{1}{|\omega_p|}
  \sum_{c\in\mathcal{C}(p)}
  \bm{q}_c\otimes\bm{\ell}_{pc}.
  \label{eq:node_gradient}
\end{align}
The dissipation coefficient $a_p$ is chosen such that
\begin{align}
a_p
\geq
\max_{\substack{c\in\mathcal{C}(p)\\k\in\{1,2\}}}
\rho_{\mathrm{sp}}\!\left(
\mathds{J}_{F,k}(\bm{q}_c)
-\bm{\mathcal{w}}_{p,k}\mathds{I}
\right), \qquad k\in\{1,2\},
\label{eq:nodal_spectral_radius}
\end{align}
where $\rho_{\mathrm{sp}}(\cdot)$ denotes the spectral radius of the flux Jacobian
\begin{align}
\mathds{J}_{F,k}(\bm{q})
=
\frac{\partial
	\big(\mathds{F}^(\bm{q})\bm{e}_k\big)}
{\partial\bm{q}},
\qquad \bm{e}_1=(1,0)^\top, \, \,  \bm{e}_2=(0,1)^\top.
\label{eq:directional_flux_jacobian}
\end{align}

For later use, we denote the integrated spatial residual by
\begin{align}
  \bm{H}_c
  :=
  -\sum_{p\in\mathcal{P}(c)}
  \widehat{\mathds{F}}_p \, \bm{\ell}_{pc},
  \qquad \textnormal{so that} \qquad 
  \frac{\mathrm{d}}{\mathrm{d}t}
  \left(|\omega_c|\bm{q}_c\right)
  =
  \bm{H}_c .
  \label{eq:integrated_residual}
\end{align}
The evolution of the cell volume is linked to the nodal grid velocities through the geometric conservation law (GCL):
\begin{align}
  \frac{\mathrm{d}}{\mathrm{d}t}|\omega_c|
  =
  \sum_{p\in\mathcal{P}(c)}
  \bm{\mathcal{w}}_p\cdot\bm{\ell}_{pc}.
  \label{eq:geometric_conservation_law}
\end{align}

\begin{schemeproperty}[Conservation and preservation of constant states]
\label{prop:node_fv_consistency}
On a periodic mesh, the semi-discrete ALE scheme \eqref{eq:node_fv_scheme} is globally conservative. Furthermore, if \eqref{eq:geometric_conservation_law} is satisfied, it preserves every spatially constant admissible state.
\end{schemeproperty}

\begin{proof}
Summing \eqref{eq:node_fv_scheme} over all cells, exchanging the cell and node sums, and using \eqref{eq:nodal_corner_closure}, we obtain
\[
  \frac{\mathrm{d}}{\mathrm{d}t} \sum_c|\omega_c|\bm{q}_c = -\sum_{p\in\mathcal{N}_h}\widehat{\mathds{F}}_p \sum_{c\in\mathcal{C}(p)}\bm{\ell}_{pc} =\bm{0}.
\]
For a constant state $\overline{\bm{q}}$, the gradient \eqref{eq:node_gradient} vanishes by \eqref{eq:nodal_corner_closure}. The physical flux contribution vanishes cellwise by \eqref{eq:corner_closure}, while the grid transport contribution is $-\overline{\bm{q}}\, \mathrm{d}|\omega_c|/\mathrm{d}t$ according to \eqref{eq:geometric_conservation_law}. It exactly balances the time derivative of $|\omega_c|\overline{\bm{q}}$.
\end{proof}

The divergence and curl constraints are measured on the dual cells by the following discrete operators
\begin{align}
  \mathfrak{D}_p(\bm{B}_h)
  =
  \sum_{c\in\mathcal{C}(p)}
  \bm{B}_c\cdot\bm{\ell}_{pc},
  \qquad
  \mathfrak{C}_p(\bm{j}_h)
  =
  \sum_{c\in\mathcal{C}(p)}
  \bm{j}_c\cdot\bm{t}_{pc}.
  \label{eq:discrete_involutions}
\end{align}
Their normalized counterparts are obtained by division by $|\omega_p|$. Up to the orientation convention induced by \eqref{eq:node_gradient}, $\mathfrak{D}_p/|\omega_p|$ and $\mathfrak{C}_p/|\omega_p|$ are finite volume approximations of $\nabla\cdot\bm{B}$ and $\nabla\times\bm{j}$, respectively. The purpose of the constrained mesh motion is to determine the interior grid velocities so that the appropriate quantity in \eqref{eq:discrete_involutions} vanishes on the current geometry configuration.

\subsection{Mesh motion and admissibility}
\label{sec.geom.motion}

The nodal coordinates satisfy the kinematic equation
\begin{align}
  \frac{\mathrm{d}\bm{x}_p}{\mathrm{d}t}=\bm{\mathcal{w}}_p,
  \qquad p\in\mathcal{N}_h.
  \label{eq:mesh_kinematics}
\end{align}
For the sake of clarity, we assume a fixed-boundary setting with $\bm{\mathcal{w}}_p=\bm{0}$ for $p\in\mathcal{N}_h^\partial$. Only the interior grid velocities are thus unknown, and the connectivity is kept fixed throughout the motion.

\begin{schemedefinition}[Admissible moving mesh]
\label{def:admissible_mesh}
For a vertex $(p,c)$ of cell $\omega_c$, let $p^+$ and $p^-$ be the next and previous vertices in the counter-clockwise ordering, and define
\begin{align}
  \mathds{E}_{pc}(\bm{X})
  :=
  [\bm{x}_{p^+}-\bm{x}_p,\,\bm{x}_{p^-}-\bm{x}_p].
\end{align}
A configuration $\bm{X}$ is admissible if every cell is simple, has nonzero edges, and is strictly convex with the prescribed counter-clockwise orientation, namely
\begin{align}
  |\omega_c(\bm{X})|>0,
  \qquad
  \det\mathds{E}_{pc}(\bm{X})>0,
  \qquad
  \lVert\bm{x}_{q}-\bm{x}_p\rVert_2>0,
  \quad q\in\{p^-,p^+\},
  \label{eq:mesh_admissibility}
\end{align}
and non-adjacent edges of the same polygon do not intersect.
\end{schemedefinition}

The characteristic mesh size $h(t)$ is given by the area-weighted quantity
\begin{align}
  h(t)
  =
  \left(
    \frac{
      \sum_{c=1}^{N_c}|\omega_c(t)|^2
    }{
      \sum_{c=1}^{N_c}|\omega_c(t)|
    }
  \right)^{1/2}.
  \label{eq:moving_mesh_size}
\end{align}

\section{Involution-preserving ALE finite volume method} \label{sec.method}
We now present the fully discrete involution-preserving scheme. We first introduce the Runge--Kutta ALE update, then derive the corresponding nonlinear constraint, its linearization, and the Lagrange-multiplier formulation. Finally, we prove that the discrete involution is preserved.

\subsection{Fully discrete Runge--Kutta ALE update}
\label{sec.method.rk}

Let $t^n$ and $\Delta t=t^{n+1}-t^n$ be the current time level and the current time step, respectively. The time step is subject to a classical CFL-type stability condition with $\textnormal{CFL}=1/2$. The time discretization relies on the SSP--RK2 method with $N_{\scriptscriptstyle\mathrm{RK}}=2$ stages and Butcher matrix $(\mathsf{a}_{sj})$ and weights $(\mathsf{b}_s)$ for $s,j=1, 
\ldots N_{\scriptscriptstyle\mathrm{RK}}$ (see \cite{gottlieb2001strong}). 

 Let $\bm{H}_c^{(s)}=\bm{H}_c(\bm{q}^{(s)},\bm{\mathcal{w}}^{(s)})$ denote the integrated finite volume residual \eqref{eq:integrated_residual} evaluated on the geometry $\bm{X}^{(s)}$ of stage $s$. Geometry and numerical solution are combined with exactly the same Runge--Kutta coefficients at each stage, namely
 	\begin{equation}
 		\bm{x}_p^{(s)}
 		=
 		\bm{x}_p^n
 		+
 		\Delta t\sum_{j=1}^{s-1}\mathsf{a}_{sj}\bm{\mathcal{w}}_p^{(j)}, \qquad 
 		|\omega_c^{(s)}| \bm{q}_c^{(s)}
 		=
 		|\omega_c^{n}| \bm{q}_c^n
 		+
 		\Delta t\sum_{j=1}^{s-1}\mathsf{a}_{sj}\bm{H}_c^{(j)},
 		\label{eq:rk_stage_geometry}
 	\end{equation}
 and the final update is
\begin{equation}
  \bm{x}_p^{n+1}
  =
  \bm{x}_p^n
  +
  \Delta t\sum_{s=1}^{N_{\scriptscriptstyle\mathrm{RK}}}
  \mathsf{b}_s\bm{\mathcal{w}}_p^{(s)},
  \qquad 
  |\omega_c^{n+1}| \bm{q}_c^{n+1}
  =
  |\omega_c^{n}| \bm{q}_c^n
  +
  \Delta t\sum_{s=1}^{N_{\scriptscriptstyle\mathrm{RK}}}
  \mathsf{b}_s\bm{H}_c^{(s)}.
  \label{eq:rk_final_state}
\end{equation}

\begin{remark}[Geometric Conservation Law]
The SSP-RK2 method is not enough to guarantee a fully discrete GCL because the polygonal area is quadratic in the vertex coordinates. Let $\bm{X}_0$ be the initial geometry and let $\alpha_j$ be the Runge--Kutta coefficient multiplying the $j^{\textnormal{th}}$ grid velocity in the targeted combination: $\alpha_j=\mathsf{a}_{sj}$ for a stage path and $\alpha_j=\mathsf{b}_j$ for the final path. Define successively
\begin{align}
  \Delta\bm{X}_j
  :=\Delta t\,\alpha_j\bm{\mathcal{w}}^{(j)},
  \qquad
  \bm{X}_j:=\bm{X}_{j-1}+\Delta\bm{X}_j.
  \label{eq:gcl_path_geometry}
\end{align}
Each nonzero displacement $\Delta\bm{X}_j$ is interpreted as a straight segment, and its geometric flux must be evaluated at the segment midpoint
\begin{align}
  \bm{X}_{j-1/2}=\bm{X}_{j-1}+\frac12\Delta\bm{X}_j,
  \qquad
  G_{c,j}^{\mathrm{mid}}
  =\sum_{p\in\mathcal{P}(c)}
  \bm{\mathcal{w}}_p^{(j)}\cdot
  \bm{\ell}_{pc}(\bm{X}_{j-1/2}).
  \label{eq:segment_midpoint_gcl_flux}
\end{align}
Since the corner vector is affine along a straight mesh trajectory, the midpoint quadrature is exact and gives
\begin{align}
  |\omega_c(\bm{X}_{j})|-|\omega_c(\bm{X}_{j-1})|
  =\Delta t\,\alpha_j G_{c,j}^{\mathrm{mid}}.
  \label{eq:segment_exact_gcl}
\end{align}
Summing \eqref{eq:segment_exact_gcl} along the actual Runge--Kutta path yields the fully discrete GCL. Cell volumes are then always recomputed from the resulting coordinates. Note that $G_{c,j}^{\mathrm{mid}}$ is a swept-geometry flux and does not necessarily coincide with the physical residual evaluated on the Runge--Kutta stage geometry.
\end{remark}

\subsection{Nonlinear constraint: involutions through mesh displacement}
\label{sec.method.target}
To provide a general description of the new methodology, let us define
\begin{align}
  \bm{z}
  =
  \begin{cases}
    \bm{B} & \textnormal{ for } \mathrm{MHD}\\
    \bm{j} & \textnormal{ for } \mathrm{EH}
  \end{cases},
  \qquad
  \bm{g}_{pc}
  =
  \begin{cases}
    \bm{\ell}_{pc} & \textnormal{ for } \mathrm{MHD}\\
    \bm{t}_{pc} & \textnormal{ for } \mathrm{EH}
  \end{cases}.
  \label{eq:constrained_variable}
\end{align}
Here $\bm{z}\in\mathbb{R}^2$ is the constrained physical field, whereas $\bm{g}_{pc}\in\mathbb{R}^2$ is the geometrical vector defining its discrete involution. In addition, $\mathds{F}_z$ denotes the contributions of the physical flux tensor $\mathds{F}$  associated with $\bm{z}$, and $\mathds{G}_{z,p}$ denotes the corresponding nodal gradient $\mathds{G}_p$. Therefore, both discrete involutions $\mathfrak{D}$ and $\mathfrak{C}$ in \eqref{eq:discrete_involutions} can be written as
\begin{align}
  \mathfrak{L}_p(\bm{z}_h,\bm{X})
  =
  \sum_{c\in\mathcal{C}(p)}
  \bm{z}_c\cdot\bm{g}_{pc} (\bm{X}), \qquad p\in\mathcal{N}_h^\circ.
  \label{eq:generic_involution}
\end{align}
The grid velocity of stage $s$ is chosen to constrain the actual Runge--Kutta combination to which this velocity contributes. To express this point, define
\begin{align}
  \beta_j^{(s)}
  =
  \begin{cases}
    \mathsf{a}_{s+1,j} & s<N_{\scriptscriptstyle\mathrm{RK}}\\
    \mathsf{b}_j & s=N_{\scriptscriptstyle\mathrm{RK}}
  \end{cases},
  \qquad
  \sigma_s:=\beta_s^{(s)}.
  \label{eq:rk_target_coefficients}
\end{align}
The sequential construction is restricted to Runge--Kutta schemes satisfying $\sigma_s>0$ for every $s$. For an intermediate stage, the target is therefore the geometry and state of the next Runge--Kutta stage; for the last stage, it is the final Runge--Kutta combination.

\begin{schemedefinition}[Runge--Kutta target of stage $s$]
\label{def:rk_target}
The contributions already known before solving for $\bm{\mathcal{w}}^{(s)}$ define the target base geometry and solution, that is
\begin{equation}
  \bm{x}_{p,\mathrm{base}}^{(s)}
  =
  \bm{x}_p^n+
  \Delta t\sum_{j=1}^{s-1}
  \beta_j^{(s)}\bm{\mathcal{w}}_p^{(j)},
  \qquad 
  |\omega_c^{(s)}| \bm{z}_{c,\mathrm{base}}^{(s)}
  =
  |\omega_c^n|\bm{z}_c^{n}
  +
  \Delta t\sum_{j=1}^{s-1}
  \beta_j^{(s)}\bm{H}_{z,c}^{(j)}.
  \label{eq:rk_target_base}
\end{equation}
The displacement unknown and the corresponding stage grid velocity are
\begin{align}
  \bm{d}_p^{(s)}
  =
  \sigma_s\Delta t\,\bm{\mathcal{w}}_p^{(s)},
  \qquad
  \bm{x}_p^{(s)}(\bm{d})
  =
  \bm{x}_{p,\mathrm{base}}^{(s)}+\bm{d}_p, \qquad p\in\mathcal{N}_h,
  \label{eq:rk_target_displacement}
\end{align}
with $\bm{d}_p=\bm{0}$ for $p\in\mathcal{N}_h^{\partial}$.
The corresponding global vectors are denoted by $\bm{d}:=(\bm{d}_p)_{p\in\mathcal{N}_h}$ and $\bm{X}^{(s)}(\bm{d}):=(\bm{x}_p^{(s)}(\bm{d}))_{p\in\mathcal{N}_h}$.
\end{schemedefinition}

The flux of stage $s$ is always evaluated on the stage geometry $\bm{X}^{(s)}$ in \eqref{eq:rk_stage_geometry}, thus in \eqref{eq:rk_target_base} we have
\begin{align}
  \bm{H}_{z,c}^{(s)}(\bm{d})
  =-
  \sum_{p\in\mathcal{P}(c)}
  \widehat{\mathds{F}}_{z,p}^{(s)}
  \left(
    \frac{\bm{d}}{\sigma_s\Delta t}
  \right)
  \bm{\ell}_{pc}^{(s)}.
  \label{eq:constrained_flux_balance}
\end{align}
The stage state, its discrete gradient, the stage corner vectors, and the numerical viscosity coefficient are fixed during the nonlinear solver. However, the ALE transport term is not frozen:
\begin{align}
  \widehat{\mathds{F}}_{z,p}^{(s)}(\bm{d})
  =
  \mathds{F}_{z}(\bm{q}_p^{(s)})
  -
  \bm{z}_p^{(s)}
  \otimes
  \frac{\bm{d}_p}{\sigma_s\Delta t}
  -
  \frac{1}{2}h_p^{(s)}\,
  a_p^{(s)}
  \mathds{G}_{z,p}^{(s)},
  \label{eq:constrained_ale_flux}
\end{align}
where the nodal value $\bm{z}_p^{(s)}$ is obtained from the stage cell averages by the volume-weighted transfer \eqref{eq:nodal_state}. It follows that the target numerical solution is
\begin{equation}
\label{eq:trial_constrained_state}
  \bm{Z}_c^{(s)}(\bm{d})
  =
  |\omega_{c,\mathrm{base}}^{(s)}| \bm{z}_{c,\mathrm{base}}^{(s)}
  +
  \sigma_s\Delta t\,\bm{H}_{z,c}^{(s)}(\bm{d}),
\end{equation}
with the definition
\begin{equation}
	\bm{Z}_c^{(j)}(\bm{d}):=|\omega_c^{(j)}(\bm{d})| \bm{z}_c^{(j)}(\bm{d}), \qquad j \in [1,N_{\scriptscriptstyle\mathrm{RK}}].
	\label{eqn.Z_def}
\end{equation}
Notice that \eqref{eq:trial_constrained_state} is recomputed at every nonlinear iteration. Freezing it would neglect the dependence of the ALE flux on the unknown grid velocity. Notice also that we directly solve for the quantity $\bm{Z}$ and not for $\bm{z}$, which can be ultimately retrieved from the definition \eqref{eqn.Z_def}.

\begin{schemedefinition}[Fully coupled involution constraint]
\label{def:coupled_constraint}
For every $p\in\mathcal{N}_h^\circ$, the displacement of stage $s$ is determined by the nonlinear equations
\begin{align}
  \mathcal{R}_p^{(s)}(\bm{d})
  :=
  \sum_{c\in\mathcal{C}(p)}
  \bm{z}_c^{(s)}(\bm{d})
  \cdot
  \bm{g}_{pc}
  \left(\bm{X}^{(s)}(\bm{d})\right)
  =0.
  \label{eq:fully_coupled_constraint}
\end{align}
Both the solution in \eqref{eq:trial_constrained_state} and all geometrical quantities entering \eqref{eq:fully_coupled_constraint} are evaluated on the current trial displacement, and solved with a Newton method.
\end{schemedefinition}


\subsection{Analytical linearization}
\label{sec.method.linearization}
Let $N_I$ denote the number of internal nodes, which are the free vertices after the boundary conditions have been imposed. Their $2N_I$ scalar displacement degrees of freedom are collected in
$\bm{d} = \left(\bm{d}_x,\bm{d}_y\right) \in\mathbb{R}^{2N_I}$. Thus $\bm{d}_x,\bm{d}_y\in\mathbb{R}^{N_I}$ are vectors, while $d_{k,r}$ is the scalar component in direction $k\in\{x,y\}$ associated with $r\in\mathcal{N}_h^\circ$. At each nonlinear iteration $m$ of the Newton solver, the residual vector and its Jacobian are denoted by
\begin{align}
  \bm{\mathcal{R}}^{(m)}
  =
  \bm{\mathcal{R}}(\bm{d}^{(m)}),
  \qquad
  \mathds{J}^{(m)}
  =
  \frac{\partial\bm{\mathcal{R}}}{\partial\bm{d}}
  \left(\bm{d}^{(m)}\right)
  =
  \begin{pmatrix}
    \mathds{J}_x^{(m)}&
    \mathds{J}_y^{(m)}
  \end{pmatrix}.
  \label{eq:constraint_jacobian}
\end{align}
There are therefore $2N_I$ scalar unknowns and, before removal of dependencies, $N_I$ scalar nodal constraints. If the constraint Jacobian has row rank $r\leq N_I$, the linearized feasible space has dimension $2N_I-r$; in the full-row-rank case its expected dimension is $N_I$. Periodic identification, prescribed boundary motion, and global patch identities may reduce $r$, whereas incompatible boundary data may make the affine linearized equations inconsistent. In practice, redundant constraint rows must be removed or handled by a rank-revealing factorization. The smallest retained singular value of the scaled Jacobian should be monitored as the mesh deforms, since its decay signals loss of regularity before cell inversion, i.e. negative volume, necessarily occurs.
The Jacobian contains three coupled effects: the variation of the corner vector, the variation of the trial cell volume, and the variation of the ALE flux. Indeed, omitting the stage index for readability, for $k\in\{x,y\}$ one has
\begin{subequations}
\label{eq:constraint_differential}
\begin{align}
  \frac{\partial\bm{z}_c}{\partial d_{k,r}}
  &=
  \frac{1}{|\omega_c|}
  \left(
    \frac{\partial\bm{Z}_c}{\partial d_{k,r}}
    -
    \bm{z}_c
    \frac{\partial|\omega_c|}{\partial d_{k,r}}
  \right),
  \label{eq:state_differential}\\
  \frac{\partial\bm{Z}_c}{\partial d_{k,r}}
  &=
  \sigma_s\Delta t\,
  \frac{\partial\bm{H}_{z,c}}{\partial d_{k,r}},
  \label{eq:numerator_differential}\\
  \frac{\partial\mathcal{R}_p}{\partial d_{k,r}}
  &=
  \sum_{c\in\mathcal{C}(p)}
  \left[
    \frac{\partial\bm{z}_c}{\partial d_{k,r}}
    \cdot\bm{g}_{pc}
    +
    \bm{z}_c\cdot
    \frac{\partial\bm{g}_{pc}}{\partial d_{k,r}}
  \right].
  \label{eq:residual_differential}
\end{align}
\end{subequations}
For a polygonal cell with counter-clockwise local node numbering, the corner vector has the geometrical representation
\begin{align}
  \bm{\ell}_{pc}
  =
  \frac{1}{2}
  \begin{pmatrix}
    y_{p^+}-y_{p^-}\\
    x_{p^-}-x_{p^+}
  \end{pmatrix},
  \qquad
  \frac{\partial|\omega_c|}{\partial\bm{x}_p}
  =
  \bm{\ell}_{pc}.
  \label{eq:geometric_derivatives}
\end{align}
Therefore its derivatives are local and nonzero only with respect to the two neighbouring vertices. The derivatives of $\bm{t}_{pc}$ follow by applying the same fixed counter-clockwise rotation. Since $\bm{q}^{(s)}$, $\mathds{G}^{(s)}$, and $a_p^{(s)}$ are frozen in \eqref{eq:constrained_ale_flux}, only the grid-transport part depends explicitly on $\bm{d}$ in the Newton method. Thus
\begin{align}
  \frac{\partial\bm{\mathcal{w}}^{(s)}_q}{\partial d_{x,r}}
  =
  \frac{\delta_{qr}}{\sigma_s\Delta t}\bm{e}_x,
  \qquad
  \frac{\partial\bm{\mathcal{w}}^{(s)}_q}{\partial d_{y,r}}
  =
  \frac{\delta_{qr}}{\sigma_s\Delta t}\bm{e}_y.
  \label{eq:velocity_displacement_derivative}
\end{align}
Here $\delta_{qr}$ is the Kronecker symbol, and $\bm{e}_x,\bm{e}_y$ are the Cartesian basis vectors. The factor $(\sigma_s\Delta t)^{-1}$ in \eqref{eq:velocity_displacement_derivative} cancels the factor $\sigma_s\Delta t$ in \eqref{eq:numerator_differential}. The resulting contribution can be written explicitly as
\begin{align}
  \left.
  \frac{\partial\bm{H}_{z,c}^{(s)}}
  {\partial d_{k,r}}
  \right|_{\mathrm{ALE}}
  =
  \frac{
    \bm{z}_r^{(s)}
    \big(\bm{\ell}_{rc}^{(s)}\big)_k
  }{\sigma_s\Delta t},
  \qquad
  k\in\{x,y\},
  \quad r\in\mathcal{P}(c).
  \label{eq:local_ale_flux_derivative}
\end{align}
It is therefore local to the cells adjacent to the displaced node, and the resulting constraint Jacobian has thus a compact node-based stencil. 

\subsection{Mesh optimization}
\label{sec.method.kkt}
The nonlinear involution constraint determines which mesh displacements are feasible, but it does not determine a unique displacement. We remove this indeterminacy with a secondary elliptic criterion. Its first part distributes the motion through the whole mesh and attracts the nodes toward the initial configuration; its second part penalizes local changes of shape and volume. Importantly, this criterion only selects the mesh: the involution itself is imposed by the fully coupled ALE constraint \eqref{eq:fully_coupled_constraint}, and no correction is applied to the physical field after the update. The construction combines linear-elastic mesh motion with the target-matrix optimization paradigm (TMOP) \cite{Tonon2021,KnuppRezoning,DobrevTMOP2019}.

\begin{schemedefinition}[Corner distortion measure]
\label{def:corner_distortion}
For every corner $(p,c)$ of cell $\omega_c$, let $p^+$ and $p^-$ be the next and previous vertices in the counter-clockwise ordering, and write $\bm{x}_p(\bm{d})=\bm{x}_{p,\mathrm{base}}^{(s)}+\bm{d}_p$. The reference and trial corner matrices are
\begin{align}
  \mathds{W}_{pc}
  &:=[\bm{x}_{p^+}^{\mathrm{ref}}-\bm{x}_p^{\mathrm{ref}},\,\bm{x}_{p^-}^{\mathrm{ref}}-\bm{x}_p^{\mathrm{ref}}],\\
  \mathds{A}_{pc}(\bm{d})
  &:=[\bm{x}_{p^+}(\bm{d})-\bm{x}_p(\bm{d}),\,\bm{x}_{p^-}(\bm{d})-\bm{x}_p(\bm{d})].
  \label{eq:corner_matrices}
\end{align}
The reference mesh is assumed admissible, so that $\det\mathds{W}_{pc}>0$ at every corner. Their relative deformation and its Jacobian are
\begin{align}
  \mathds{T}_{pc}(\bm{d})
  :=\mathds{A}_{pc}(\bm{d})\mathds{W}_{pc}^{-1},
  \qquad
  J_{pc}(\bm{d}):=\det\mathds{T}_{pc}(\bm{d}).
  \label{eq:corner_deformation_tensor}
\end{align}
For $J=\det\mathds{T}>0$, define
\begin{align}
  \varrho_{\mathrm{q}}(\mathds{T})
  &:=
  \frac{(T_{11}-T_{22})^2+(T_{12}+T_{21})^2}{2J}
  +\frac{3}{4}\left(\sqrt{J}-\frac{1}{\sqrt{J}}\right)^2,
  &
  \Psi(\mathds{T})&:=\varrho_{\mathrm{q}}(\mathds{T})+\varrho_{\mathrm{q}}(\mathds{T})^2,
  \label{eq:corner_quality_density}
\end{align}
and
\begin{align}
  \mathcal{E}_{\mathrm{q}}(\bm{d})
  :=\sum_{c=1}^{N_c}\frac{1}{|\mathcal{P}(c)|}
  \sum_{p\in\mathcal{P}(c)}\det(\mathds{W}_{pc})\,
  \Psi\!\left(\mathds{T}_{pc}(\bm{d})\right).
  \label{eq:mesh_quality_energy}
\end{align}
The first contribution in \eqref{eq:corner_quality_density} measures shear and anisotropy, while the second measures the local change of area. Hence $\Psi$ in \eqref{eq:mesh_quality_energy} vanishes at the reference configuration and diverges when a corner approaches degeneracy.
\end{schemedefinition}

Let $\bm{X}^{\mathrm{ref}}$ be the initial mesh and let $a_{\max}^{(s)}$ be the largest stage wave-speed bound. We introduce the dimensionless stage Courant number and the distortion of the base mesh, that is
\begin{align}
  \chi_s
  &:=
  \frac{\sigma_s\Delta t\,a_{\max}^{(s)}}{h_{\min}^{\mathrm{base}}},
  &
  h_{\min}^{\mathrm{base}}
  &:=
  \min_c\sqrt{|\omega_c^{\mathrm{base}}|},
  \\
  \delta_s
  &:=
  \max_{(p,c)}\sqrt{2\Psi\!\left(\mathds{T}_{pc}(\bm{0})\right)},
  \label{eq:predictor_measures}
\end{align}
which measures the largest corner distortion of the base mesh. We further define the bounded function
\begin{align}
  \zeta(\chi)
  :=
  \begin{cases}
    \dfrac{\chi}{1-\chi}, & 0\leq\chi<\dfrac12\\[1ex]
    1, & \chi\geq\dfrac12
  \end{cases},
\end{align}
and set
\begin{align}
  \bm{d}_{p}^{\mathrm{ref},(s)} :=\vartheta_s\left(\bm{x}_{p}^{\mathrm{ref}}-\bm{x}_{p,\mathrm{base}}^{(s)}\right),
  \qquad
  \bm{e}_p:=\bm{d}_p-\bm{d}_{p}^{\mathrm{ref},(s)},
  \label{eq:reference_displacement}
\end{align}
with the definition 
\begin{eqnarray}
	\vartheta_s
	:=
	\frac{\chi_s}{1+\chi_s}
	+\frac{\zeta(\chi_s)}{1+\chi_s}
	\frac{\delta_s}{\chi_s+\delta_s},
\end{eqnarray}
and we take $\vartheta_s=0$ when $\chi_s=\delta_s=0$. The first term in $\vartheta_s$ supplies a bounded return toward the reference mesh on the stage time scale; the second term activates a stronger recovery only when distortion has accumulated. All quantities are computed from the current stage. For a fixed mesh, bounded wave speed, and bounded distance between the base and reference coordinates, $\bm{d}^{\mathrm{ref},(s)}=O(\Delta t)$ as $\Delta t\to0$.

For a nodal vector field $\bm{v}_h=(\bm{v}_p)_{p\in\mathcal{N}_h}$ satisfying $\bm{v}_p=\bm{0}$ for $p\in\mathcal{N}_h^\partial$, define the cellwise gradient and symmetric strain
\begin{align}
  \nabla_c^h\bm{v}_h
  &:=
  \frac{1}{|\omega_c^{\mathrm{base}}|}
  \sum_{r\in\mathcal{P}(c)}\bm{v}_r\otimes\bm{\ell}_{rc}^{\mathrm{base}},
  &
  \bm{\varepsilon}_c(\bm{v}_h)
  &:=
  \frac{1}{2}\left(\nabla_c^h\bm{v}_h+(\nabla_c^h\bm{v}_h)^\top\right).
  \label{eq:discrete_mesh_strain}
\end{align}
With $L_\Omega$ the largest Cartesian extent of the domain, the elliptic bilinear form is
\begin{align}
  \mathcal{a}_{\mathrm{el}}(\bm{u}_h,\bm{v}_h)
  :=\sum_{c=1}^{N_c}\left[
  \frac{|\omega_c^{\mathrm{base}}|}{|\mathcal{P}(c)|}
  \sum_{p\in\mathcal{P}(c)}\bm{u}_p\cdot\bm{v}_p
  +2L_\Omega^2|\omega_c^{\mathrm{base}}|\,
  \bm{\varepsilon}_c(\bm{u}_h):\bm{\varepsilon}_c(\bm{v}_h)\right].
  \label{eq:anchored_elastic_form}
\end{align}
It discretizes the screened elasticity operator
\begin{align}
  \bm{e}-L_\Omega^2\nabla\cdot\left(2\bm{\varepsilon}(\bm{e})\right), \qquad \bm{e} = \bm{d} - \bm{d}^{\mathrm{ref}}.
  \label{eq:anchored_elastic_operator}
\end{align}
The zeroth-order term removes rigid modes and anchors the displacement, whereas the strain term couples its two components and distributes a local correction over neighbouring cells. The factor $L_\Omega^2$ normalizes the elastic coupling with the domain scale; no mesh-dependent tuning parameter is here introduced.

\begin{schemedefinition}[Constrained mesh problem]
\label{def:constrained_mesh_problem}
At stage $s$, the mesh-selection problem is
\begin{subequations}
\label{eq:constrained_mesh_problem}
\begin{align}
  \min_{\bm{d}}\quad
  \mathcal{E}_{h}^{(s)}(\bm{d})
  &:=
  \frac{1}{2}\mathcal{a}_{\mathrm{el}}(\bm{e}_h,\bm{e}_h)
  +\chi_s\mathcal{E}_{\mathrm{q}}(\bm{d}),
  \label{eq:constrained_mesh_objective}\\
  \text{subject to}\quad
  \bm{\mathcal{R}}^{(s)}(\bm{d})&=\bm{0},
  \qquad \text{with} \quad 
  \bm{X}^{(s)}(\bm{d})\ \text{admissible}.
\end{align}
\end{subequations}
\end{schemedefinition}

A feasible Sequential Quadratic Programming (SQP) correction is obtained as follows. At iteration $m$, let $\mathds{H}^{(m)}$ be the Gauss--Newton model of \eqref{eq:constrained_mesh_objective}, made positive definite by the elliptic part. We denote the objective gradient and the Jacobian of the involution constraint by
\begin{align}
  \bm{\gamma}^{(m)}
  &:=(\bm{\gamma}_x^{(m)},\bm{\gamma}_y^{(m)})
  =\nabla_{\!\bm{d}}\mathcal{E}_h^{(s)}(\bm{d}^{(m)}),
  &
  \mathds{J}^{(m)}
  &:=(\mathds{J}_x^{(m)},\mathds{J}_y^{(m)}).
\end{align}
The displacement correction and Lagrange multiplier solve
\begin{align}
  \begin{pmatrix}
    \mathds{H}_{xx}&\mathds{H}_{xy}&-\mathds{J}_x^\top\\
    \mathds{H}_{yx}&\mathds{H}_{yy}&-\mathds{J}_y^\top\\
    \mathds{J}_x&\mathds{J}_y&\bm{0}
  \end{pmatrix}
  \begin{pmatrix}
    \delta\bm{d}_x\\
    \delta\bm{d}_y\\
    \bm{\lambda}
  \end{pmatrix}
  =-
  \begin{pmatrix}
    \bm{\gamma}_x\\
    \bm{\gamma}_y\\
    \bm{\mathcal{R}}
  \end{pmatrix}.
  \label{eq:kkt_system}
\end{align}
All blocks and right-hand sides in \eqref{eq:kkt_system} are evaluated at iteration $m$; the superscript is omitted there for readability. The vector $\bm{\lambda}\in\mathbb{R}^{N_I}$ contains the Lagrange multipliers associated with the nodal involution constraints. Note that the elliptic metric and the constraint Jacobian have compact node-based stencils; consequently, \eqref{eq:kkt_system} is a sparse saddle-point system.
During feasibility restoration, $\bm{\gamma}_x=\bm{\gamma}_y=\bm{0}$; after feasibility has been reached, the last right-hand side is set to zero and the same system produces a tangent mesh-quality step.

\begin{schemeproposition}[Elliptic selection of a feasible correction]
\label{prop:elliptic_selection}
Assume that $\mathds{H}^{(m)}$ is positive definite on the free displacement space and that $\mathds{J}^{(m)}$ has full row rank. Then \eqref{eq:kkt_system} has a unique displacement component. During restoration, this component is the minimum-$\mathds{H}^{(m)}$-norm correction satisfying the linearized involution constraint.
\end{schemeproposition}

\begin{proof}
The first two block rows are the stationarity equations of a strictly convex quadratic functional under the affine constraint given by the last block row. The KKT conditions are therefore necessary and sufficient, and strict convexity gives a unique displacement correction.
\end{proof}

The SQP iteration is globalized by a backtracking line search restricted to admissible meshes. Restoration corrections first recover the nonlinear involution constraint; tangent corrections then decrease \eqref{eq:constrained_mesh_objective} to first order along the constraint manifold, with feasibility restored after each accepted correction. Therefore the procedure seeks an admissible local solution of \eqref{eq:constrained_mesh_problem}, while the accepted state satisfies the complete non-frozen constraint to the prescribed nonlinear tolerance.

\begin{schemedefinition}[Predictor--corrector algorithm]
An inexpensive initial displacement is first obtained by solving a frozen-flux problem, in which the numerical flux is evaluated only once while the trial geometry is allowed to vary. The subsequent fully coupled corrector updates the geometry, ALE flux, conservative state, and involution residual at each nonlinear iteration. Only solutions satisfying the complete, non-frozen nonlinear problem are considered admissible.
\end{schemedefinition}


\begin{algorithm}[H]
\caption{Involution-preserving mesh displacement at Runge--Kutta stage $s$.}
\label{alg:constrained_mesh_motion}
\begin{algorithmic}[1]
\STATE Construct the Runge--Kutta target base geometry and conservative
state from \eqref{eq:rk_target_base}.
\STATE Compute an initial displacement by solving a frozen-flux version of
\eqref{eq:constrained_flux_balance}--\eqref{eq:fully_coupled_constraint}.
\STATE Solve the constrained mesh problem \eqref{eq:constrained_mesh_problem}
by successive KKT corrections \eqref{eq:kkt_system}. At every iterate,
recompute the target geometry, ALE flux \eqref{eq:constrained_ale_flux},
conservative state \eqref{eq:trial_constrained_state}, and constraint
Jacobian \eqref{eq:constraint_jacobian}.
\STATE Globalize each correction by a line search restricted to meshes
satisfying \eqref{eq:mesh_admissibility}.
\STATE Accept the displacement when the complete non-frozen residual
\eqref{eq:fully_coupled_constraint} reaches the prescribed tolerance.
\STATE Recover the grid velocity from \eqref{eq:rk_target_displacement} and
use the same constrained flux balance in the Runge--Kutta update
\eqref{eq:rk_stage_geometry}--\eqref{eq:rk_final_state}.
\end{algorithmic}
\end{algorithm}

\subsection{Discrete involution preservation}
\label{sec.method.preservation}

\begin{schemetheorem}[Preservation of the targeted Runge--Kutta state]
\label{thm:rk_involution_preservation}
Assume that the state at time $t^n$ satisfies the appropriate discrete involution. At every Runge--Kutta stage, suppose that \eqref{eq:fully_coupled_constraint} is solved exactly and that the flux balance used in \eqref{eq:trial_constrained_state} is exactly the flux balance entering the finite volume residual. Then every targeted intermediate state and the final state \eqref{eq:rk_final_state} satisfy the discrete involution exactly. With an iterative nonlinear solver, the constraint is preserved up to its stopping tolerance and floating-point roundoff.
\end{schemetheorem}

\begin{proof}
At stage $s$, combine \eqref{eq:rk_target_base} with the current contribution $\sigma_s\Delta t\,\bm{H}_{z,c}^{(s)}$. Because the same flux balance is used by the constrained solver and the finite volume update, the resulting numerator is precisely \eqref{eq:trial_constrained_state}. Likewise, \eqref{eq:rk_target_base} and \eqref{eq:rk_target_displacement} produce exactly the geometry of the targeted Runge--Kutta combination. Therefore the cell average entering that combination is \eqref{eq:trial_constrained_state}, not an auxiliary state. Equation \eqref{eq:fully_coupled_constraint} is consequently identical to the discrete involution \eqref{eq:generic_involution} evaluated on the targeted state. For $s<N_{\scriptscriptstyle\mathrm{RK}}$ this establishes the property for the next Runge--Kutta stage. At $s=N_{\scriptscriptstyle\mathrm{RK}}$, $\beta_j^{(s)}=\mathsf{b}_j$, so the target is exactly \eqref{eq:rk_final_state}. The conclusion follows stage by stage.
\end{proof}

\begin{schemeproposition}[Inexact nonlinear solves]
\label{prop:inexact_preservation}
Suppose the accepted displacement at target $s$ satisfies
\begin{align}
  \lVert\bm{\mathcal{R}}^{(s)}\rVert_\star\leq\tau_s, \qquad \tau_s > 0,
  \label{eq:stage_nonlinear_tolerance}
\end{align}
in a fixed norm $\lVert\cdot\rVert_\star$, and that the state actually updated by the finite volume code is the same state used to evaluate this residual. Then the targeted discrete involution satisfies the same bound. In particular,
\begin{align}
  \lVert\bm{\mathfrak{L}}^{n+1}\rVert_\star
  \leq \tau_{N_{\scriptscriptstyle\mathrm{RK}}}+\varepsilon_{\mathrm{alg}},
  \label{eq:inexact_final_bound}
\end{align}
where $\varepsilon_{\mathrm{alg}}$ collects floating-point contributions. Hence one may take $C_{\mathrm{RK}}=1$ for the RK-target formulation; errors from earlier targets do not accumulate in \eqref{eq:inexact_final_bound} because the last nonlinear residual is evaluated directly on the final RK combination.
\end{schemeproposition}

For a normalized diagnostic, the same conclusion holds after division by the smallest interior nodal-patch area. Thus the raw tolerance should scale like the desired normalized accuracy times $\min_{p\in\mathcal{N}_h^\circ}|\omega_p|$ (and, in convergence studies, should be chosen below the targeted truncation error). This provides a mesh-dependent stopping rule instead of a fixed raw residual threshold.

\begin{remark}[Scope and modularity of this construction]
\label{rem:scope_modularity}
Theorem \ref{thm:rk_involution_preservation} applies equally to the divergence constraint of ideal MHD and to the curl constraint of EH. The two cases differ only through the constrained variable, the geometrical vector in \eqref{eq:constrained_variable}, and the corresponding physical flux. The Runge--Kutta target, nonlinear coupling, admissibility checks, and linear algebra are common to both. 

The proof does not rely on the particular Rusanov-type flux in \eqref{eq:stabilized_nodal_flux}. The construction can be coupled to another conservative ALE spatial operator provided that its discrete involution is available on the updated state and geometry, the identical constrained residual is used in both \eqref{eq:fully_coupled_constraint} and the conservative update, and every mesh-dependent contribution is included in the nonlinear residual and its linearization. 
Under these conditions, the constrained mesh solve is a modular structure-preserving layer rather than a model-specific field correction.
\end{remark}

\section{Numerical results} \label{sec.results}

In order to validate the scheme presented in this paper, we solve in the following a set of classical benchmark test cases. A Cartesian computational mesh is initially assigned, which then evolves in order to respect the prescribed involutions.

\subsection{MHD numerical convergence studies}

The numerical convergence studies for the new involution-preserving scheme are carried out using a smooth manufactured vortex solution defined over the computational domain $\Omega=[-5,5]^2$, starting from \cite{balsara1999staggered}. The exact solution is given by 
\begin{equation}
	\begin{aligned}
		\rho = 1 + 0.1\,\psi,
		\qquad
			\bm{u} =
		\begin{pmatrix}
			1 \\ 1
		\end{pmatrix}
		+
		\frac{\varepsilon}{2\pi}\psi
		\begin{pmatrix}
			-y\\
			x
		\end{pmatrix},
		\qquad
			\mathbf{B} =
		\frac{\mu}{2\pi}\psi
		\begin{pmatrix}
			-y\\
			x
		\end{pmatrix},
	\end{aligned}
\end{equation}
with
\begin{equation}
	\psi=\exp\left(\frac{1-r^2}{2}\right),
	\qquad
	r=\sqrt{x^2+y^2},
\end{equation}
while the pressure is defined as
\begin{equation}
	p=\varepsilon+\frac{\mu}{\sqrt{4\pi}}
	+\frac{1}{8\pi}\left(\frac{\mu}{2\pi}\right)^2
	(1-r^2)e^{1-r^2}
	-\frac{1}{2}\left(\frac{\varepsilon}{2\pi}\right)^2
	e^{1-r^2}
	+0.1\,\psi.
\end{equation}
The parameters are set to $\varepsilon=1$ and $\mu=4\pi$.

The simulations are run on a sequence of successively refined meshes until the final time $t_f = 0.1$. In Table \ref{tab:mhd_convergence}, we report the errors in $L^1$ and $L^2$ norm. We observe that the expected order of one is achieved.

\begin{table}[htbp]
	\centering
	\caption{MHD convergence rates on moving meshes at $t_f=0.1$.
		$L^1$ and $L^2$ errors and experimental orders of convergence for the primitive variables.}
	\label{tab:mhd_convergence}
	\scriptsize
	\setlength{\tabcolsep}{3pt}
	
	\begin{tabular}{c cc cc cc cc cc cc}
		\hline
		$N_x=N_y$
		& \multicolumn{2}{c}{$\rho$}
		& \multicolumn{2}{c}{$u$}
		& \multicolumn{2}{c}{$v$}
		& \multicolumn{2}{c}{$p$}
		& \multicolumn{2}{c}{$B_x$}
		& \multicolumn{2}{c}{$B_y$} \\
		\hline
		
		\multicolumn{13}{c}{$L^1$ errors} \\
		\hline
		
		30
		& $2.017\mathrm{e}{-3}$ & -
		& $3.379\mathrm{e}{-3}$ & -
		& $3.369\mathrm{e}{-3}$ & -
		& $4.863\mathrm{e}{-3}$ & -
		& $4.234\mathrm{e}{-2}$ & -
		& $4.241\mathrm{e}{-2}$ & - \\
		
		40
		& $1.532\mathrm{e}{-3}$ & 0.964
		& $2.557\mathrm{e}{-3}$ & 0.978
		& $2.564\mathrm{e}{-3}$ & 0.958
		& $3.684\mathrm{e}{-3}$ & 0.973
		& $3.206\mathrm{e}{-2}$ & 0.975
		& $3.222\mathrm{e}{-2}$ & 0.963 \\
		
		50
		& $1.232\mathrm{e}{-3}$ & 0.989
		& $2.059\mathrm{e}{-3}$ & 0.984
		& $2.071\mathrm{e}{-3}$ & 0.969
		& $2.997\mathrm{e}{-3}$ & 0.938
		& $2.580\mathrm{e}{-2}$ & 0.988
		& $2.599\mathrm{e}{-2}$ & 0.976 \\
		
		60
		& $1.058\mathrm{e}{-3}$ & 0.840
		& $1.757\mathrm{e}{-3}$ & 0.876
		& $1.768\mathrm{e}{-3}$ & 0.873
		& $2.607\mathrm{e}{-3}$ & 0.768
		& $2.169\mathrm{e}{-2}$ & 0.957
		& $2.186\mathrm{e}{-2}$ & 0.955 \\
		
		\hline
		\multicolumn{13}{c}{$L^2$ errors} \\
		\hline
		
		30
		& $5.006\mathrm{e}{-3}$ & -
		& $8.167\mathrm{e}{-3}$ & -
		& $8.123\mathrm{e}{-3}$ & -
		& $1.805\mathrm{e}{-2}$ & -
		& $1.013\mathrm{e}{-1}$ & -
		& $1.017\mathrm{e}{-1}$ & - \\
		
		40
		& $3.818\mathrm{e}{-3}$ & 0.950
		& $6.207\mathrm{e}{-3}$ & 0.962
		& $6.184\mathrm{e}{-3}$ & 0.956
		& $1.398\mathrm{e}{-2}$ & 0.895
		& $7.718\mathrm{e}{-2}$ & 0.955
		& $7.754\mathrm{e}{-2}$ & 0.951 \\
		
		50
		& $3.070\mathrm{e}{-3}$ & 0.991
		& $5.011\mathrm{e}{-3}$ & 0.973
		& $5.008\mathrm{e}{-3}$ & 0.958
		& $1.139\mathrm{e}{-2}$ & 0.932
		& $6.236\mathrm{e}{-2}$ & 0.968
		& $6.278\mathrm{e}{-2}$ & 0.959 \\
		
		60
		& $2.589\mathrm{e}{-3}$ & 0.938
		& $4.230\mathrm{e}{-3}$ & 0.934
		& $4.226\mathrm{e}{-3}$ & 0.937
		& $9.628\mathrm{e}{-3}$ & 0.926
		& $5.258\mathrm{e}{-2}$ & 0.941
		& $5.288\mathrm{e}{-2}$ & 0.947 \\
		
		\hline
	\end{tabular}
\end{table}

%

\subsection{MHD Rotor Problem}
\begin{figure}[t]
\centering
\begin{minipage}[t]{0.32\textwidth}\centering
\includegraphics[width=\linewidth]{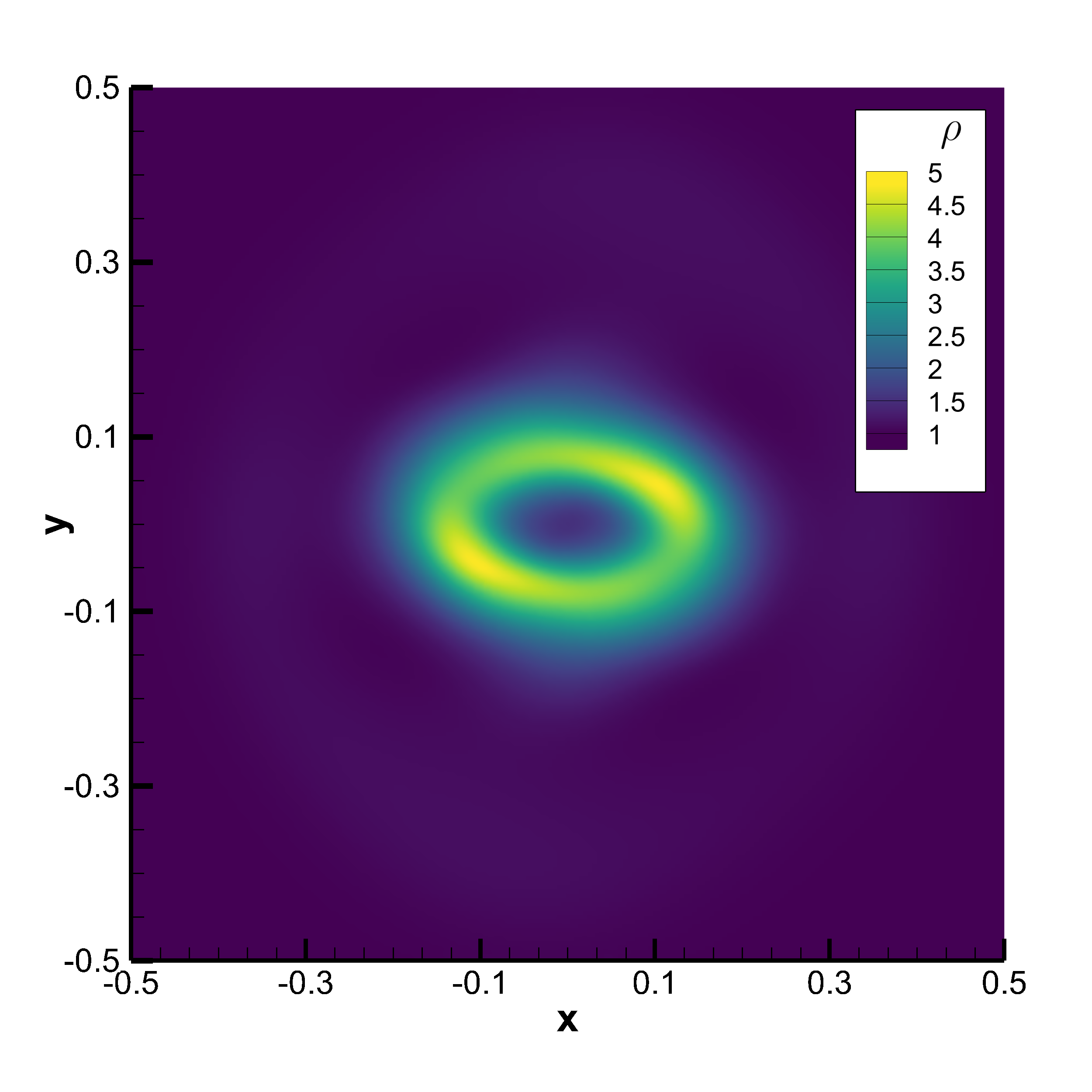}\\[-1mm]
\footnotesize{(a) Density}
\end{minipage}\hfill
\begin{minipage}[t]{0.32\textwidth}\centering
\includegraphics[width=\linewidth]{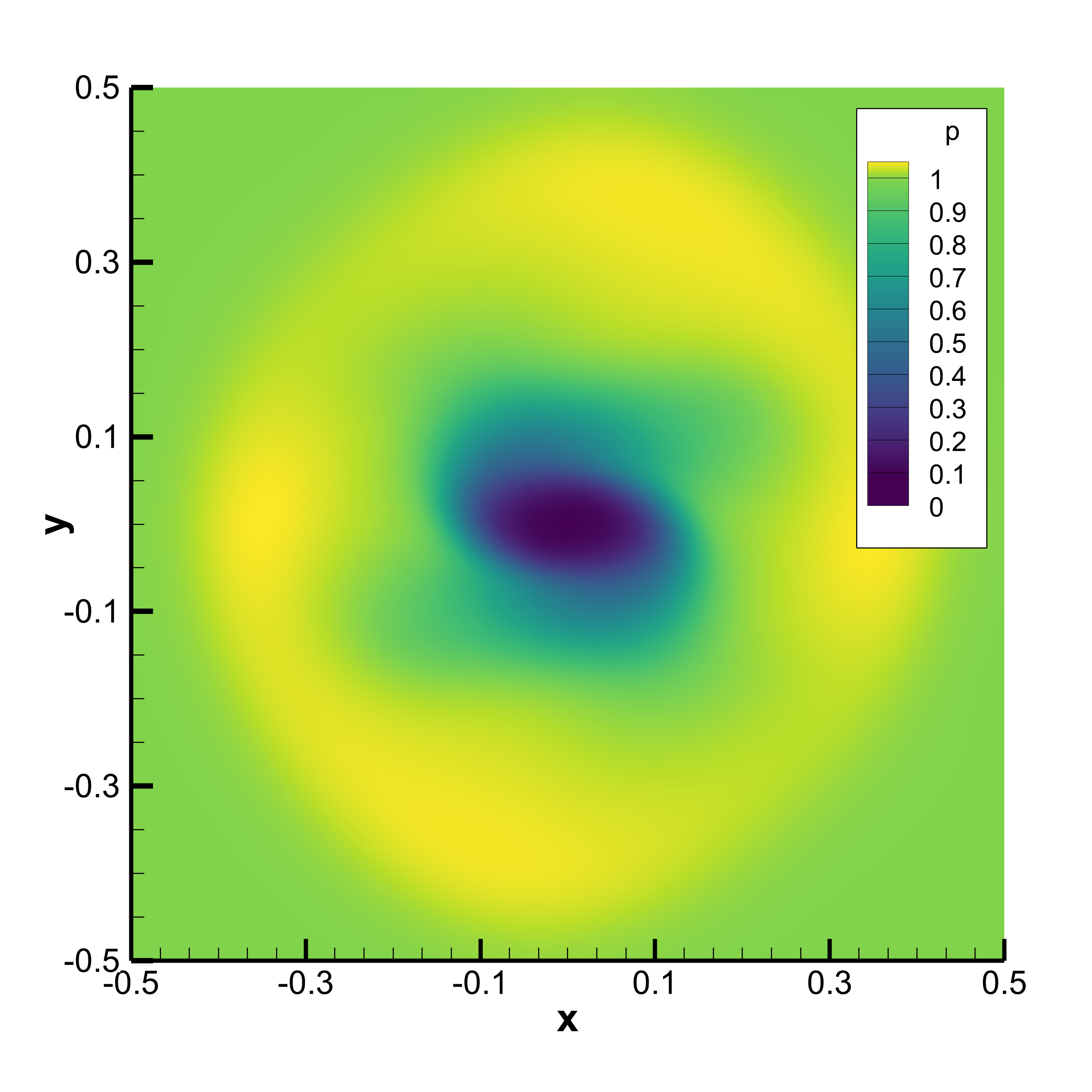}\\[-1mm]
\footnotesize{(b) Pressure}
\end{minipage}\hfill
\begin{minipage}[t]{0.32\textwidth}\centering
\includegraphics[width=\linewidth]{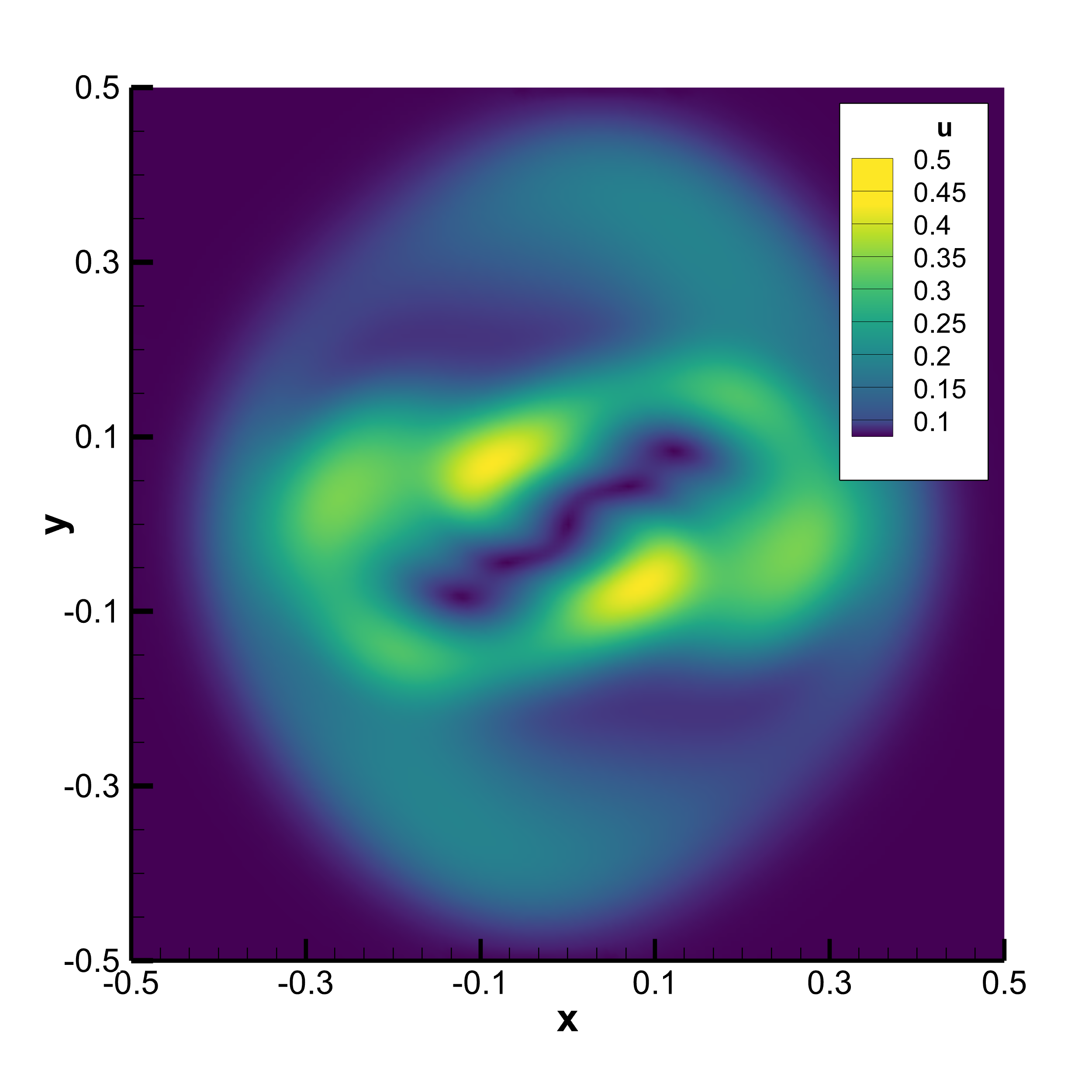}\\[-1mm]
\footnotesize{(c) Total velocity}
\end{minipage}

\medskip

\begin{minipage}[t]{0.32\textwidth}\centering
\includegraphics[width=\linewidth]{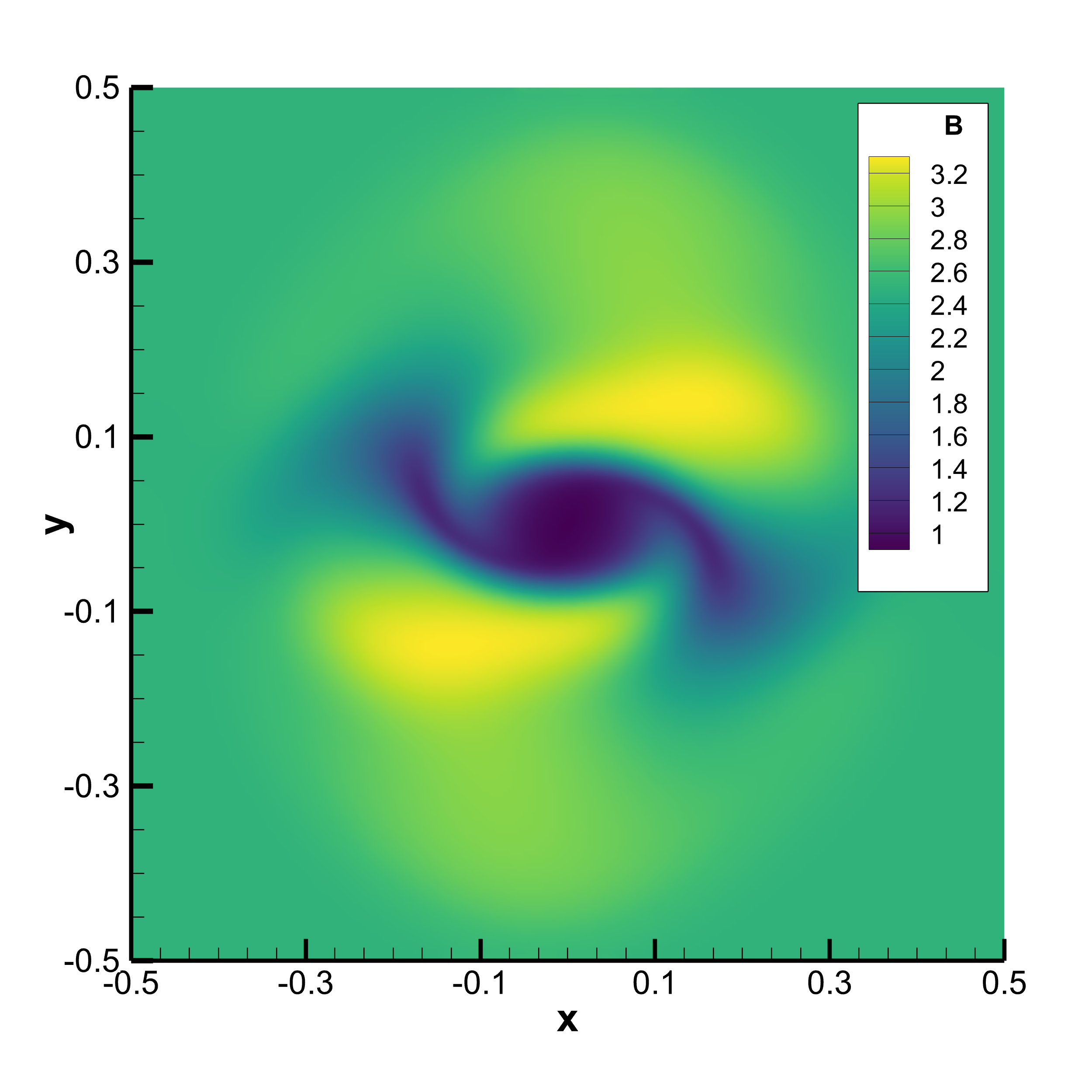}\\[-1mm]
\footnotesize{(d) Total magnetic field}
\end{minipage}
\begin{minipage}[t]{0.32\textwidth}\centering
\includegraphics[width=\linewidth]{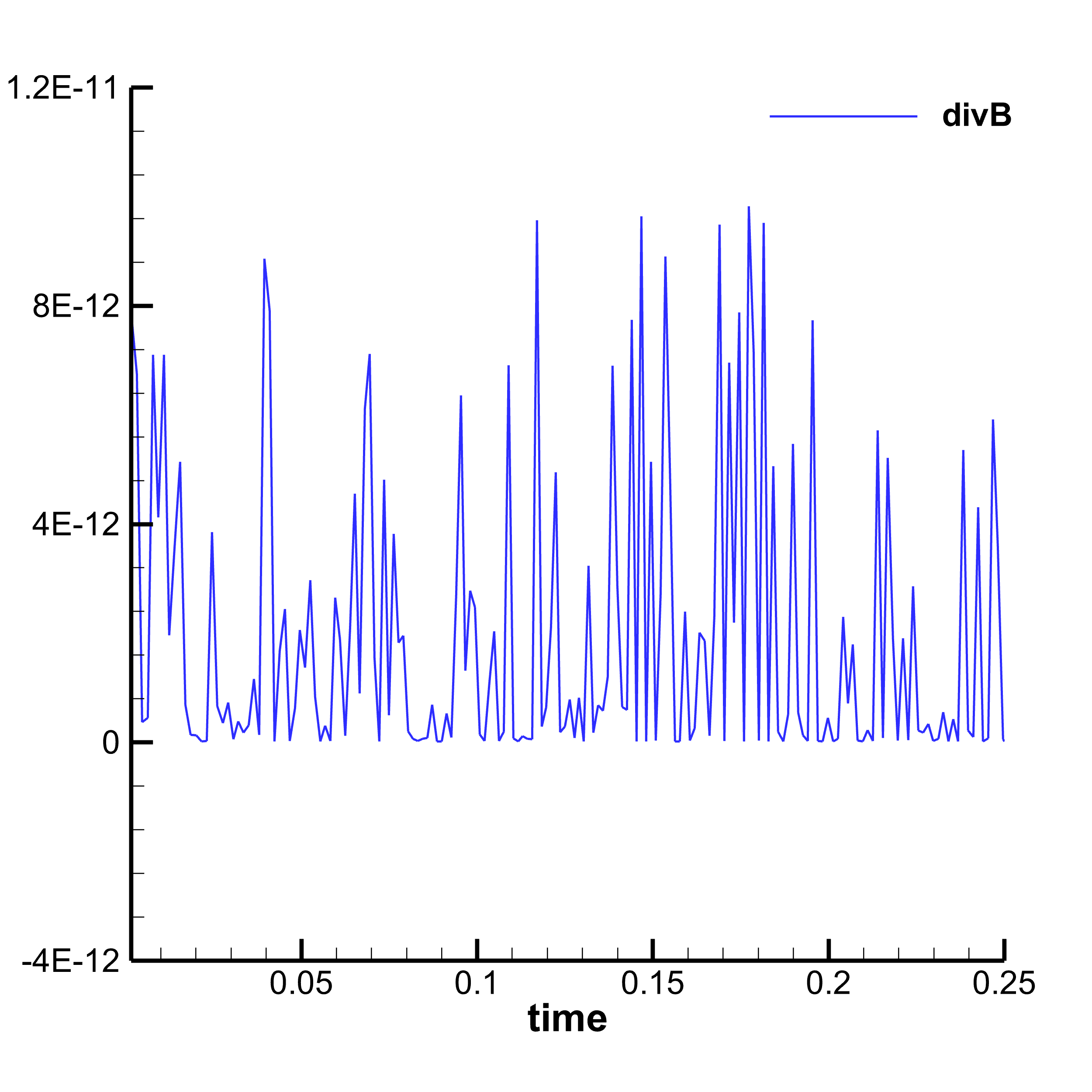}\\[-1mm]
\footnotesize{(e) Time evolution of $\nabla \cdot \mathbf{B}$}
\end{minipage}
\caption{Numerical solution of the \textit{MHD Rotor problem} at time $t_f = 0.25$ obtained using the new involution-preserving scheme on a $Nx = Ny = 200$ mesh. }
\label{fig:mhd_rotor}
\end{figure}

\begin{figure}[t]
\centering
\begin{minipage}[t]{0.48\textwidth}\centering
\includegraphics[width=0.85\linewidth]{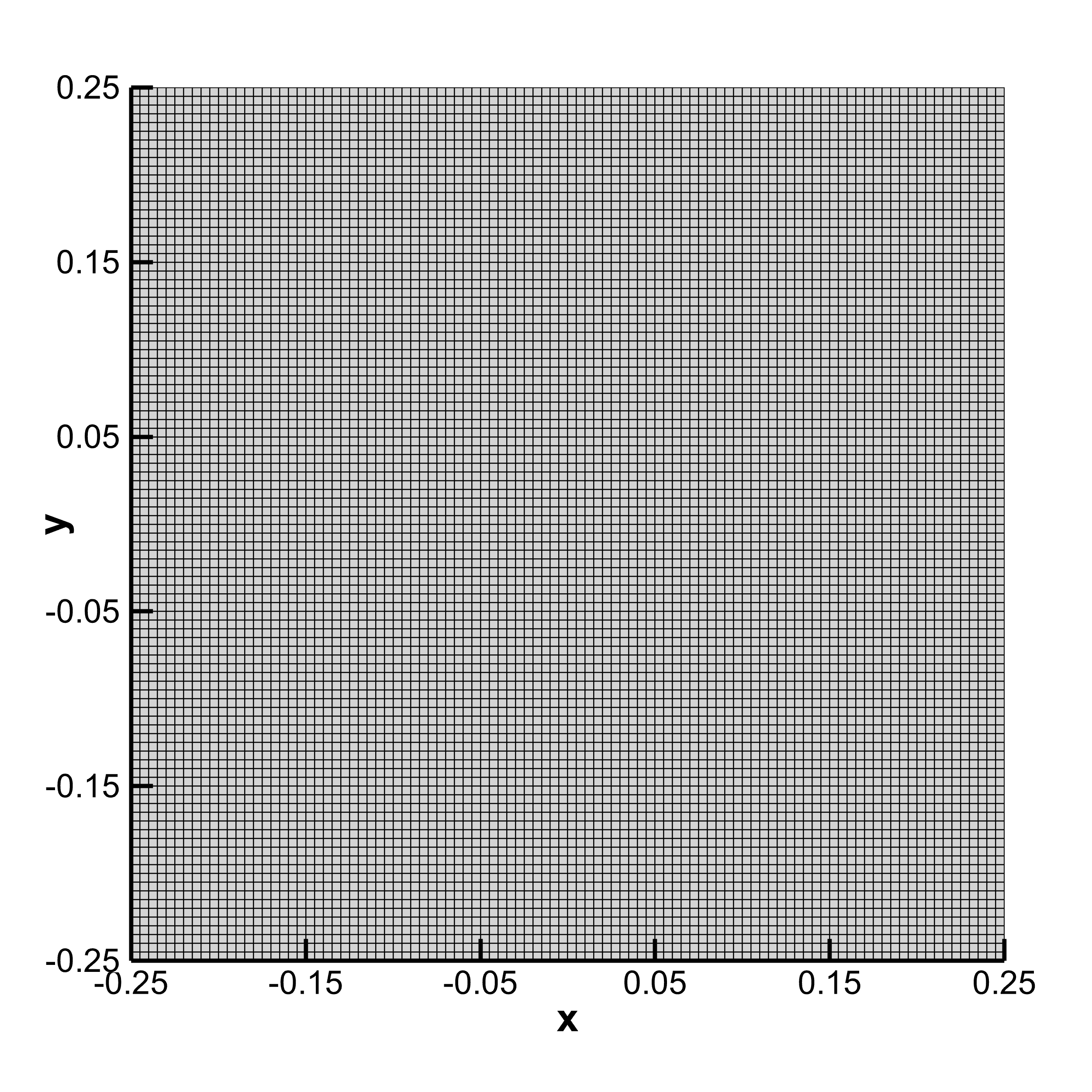}\\[-1mm]
{\footnotesize (a) Mesh at $t=0$}
\end{minipage}\hfill
\begin{minipage}[t]{0.48\textwidth}\centering
\includegraphics[width=0.85\linewidth]{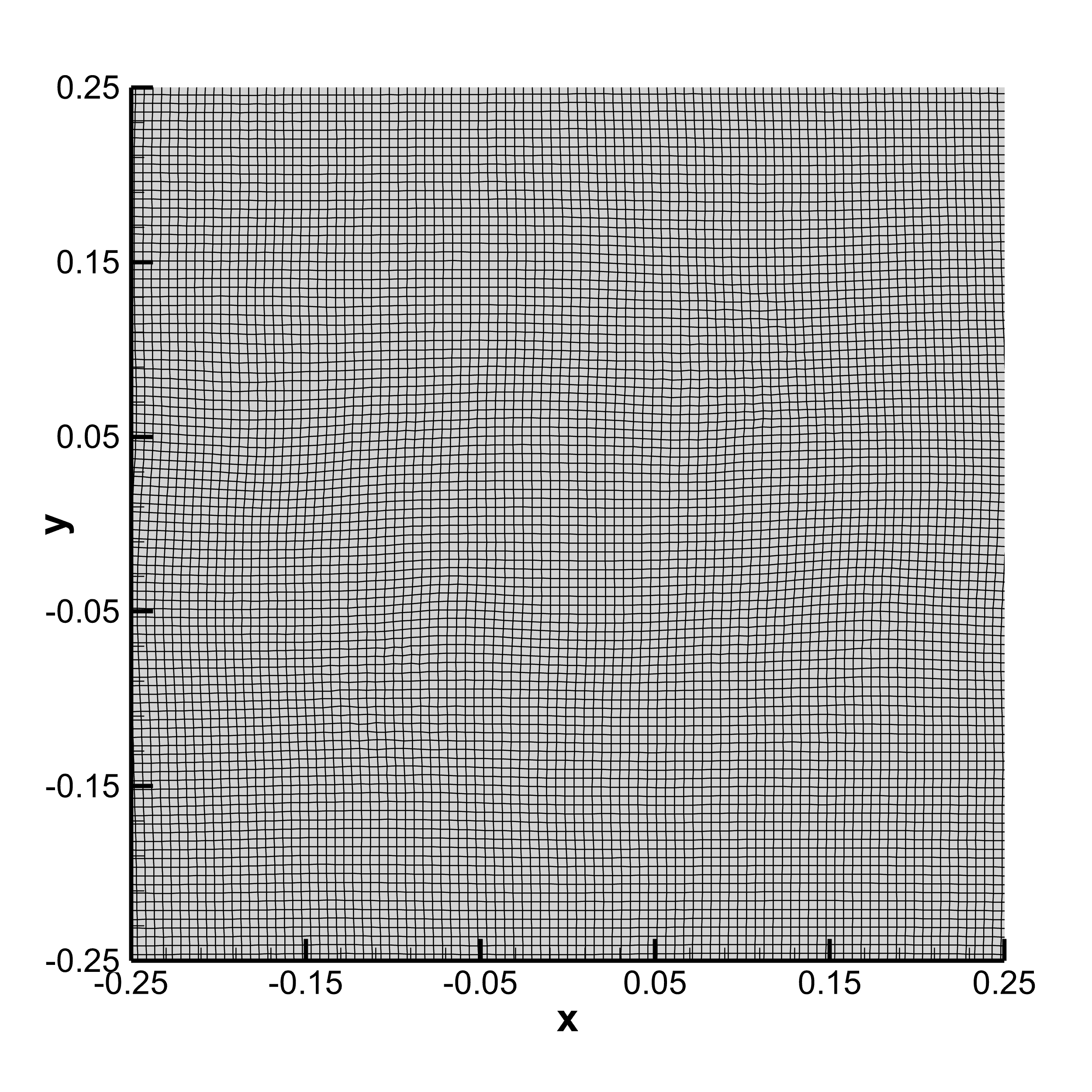}\\[-1mm]
{\footnotesize (b) Mesh at $t=t_f$}
\end{minipage}
\caption{Central zoom $[-0.25,0.25]^2$ of the rotor mesh at the initial and final time. The comparison highlights the smooth redistribution of the ALE displacement without obscuring the cell structure of the $N_x = N_y = 200$ mesh.}
\label{fig:mhd_rotor_mesh}
\end{figure}

The so-called rotor problem proposed in \cite{balsara1999staggered} is characterized by a large density difference inside and outside the rotor. The computational domain is the square $\Omega = [-0.5, 0.5]^2$, and $r:=\sqrt{x^2+y^2}$ denotes the distance from its centre. The background pressure and magnetic field are constant and set to $p = 1$ and $\bm{B} = (2.5, 0)^\top$, respectively.
The density and the velocity field are defined differently inside and outside the rotor:
\begin{equation}
	\begin{aligned}
		\rho &=
		\begin{cases}
			10 & 0 \le r \le 0.1,\\
			1  & \text{otherwise},
		\end{cases}
		\qquad
		(u,v) &=
		\begin{cases}
			(-10y,10x) & 0 \le r \le 0.1,\\
			(0,0) & \text{otherwise}.
		\end{cases}
	\end{aligned}
\end{equation}
The ratio of specific heats is set to $\gamma = 1.4$ and the simulation is run until the final time $t_f = 0.25$. The numerical solution is displayed in Figure~\ref{fig:mhd_rotor}, while Figure~\ref{fig:mhd_rotor_mesh} compares the initial and final meshes in the central region. Over the complete computation, the maximum discrete divergence is $9.82\times10^{-12}$ and its final value is $2.06\times10^{-14}$.

\subsection{MHD blast wave}
\begin{figure}[t]
\centering
\begin{minipage}[t]{0.32\textwidth}\centering
\includegraphics[width=\linewidth]{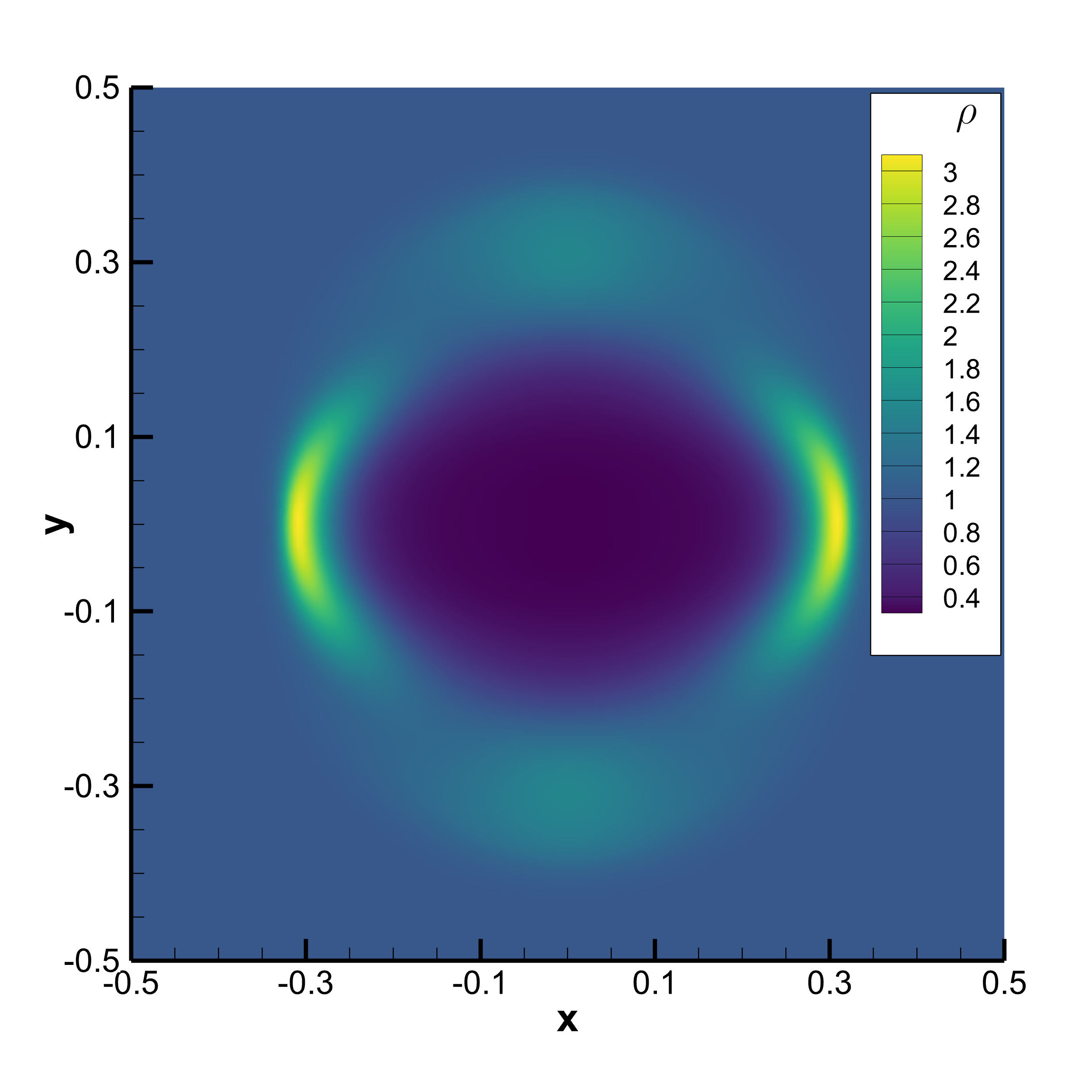}\\[-1mm]
\footnotesize{(a) Density}
\end{minipage}\hfill
\begin{minipage}[t]{0.32\textwidth}\centering
\includegraphics[width=\linewidth]{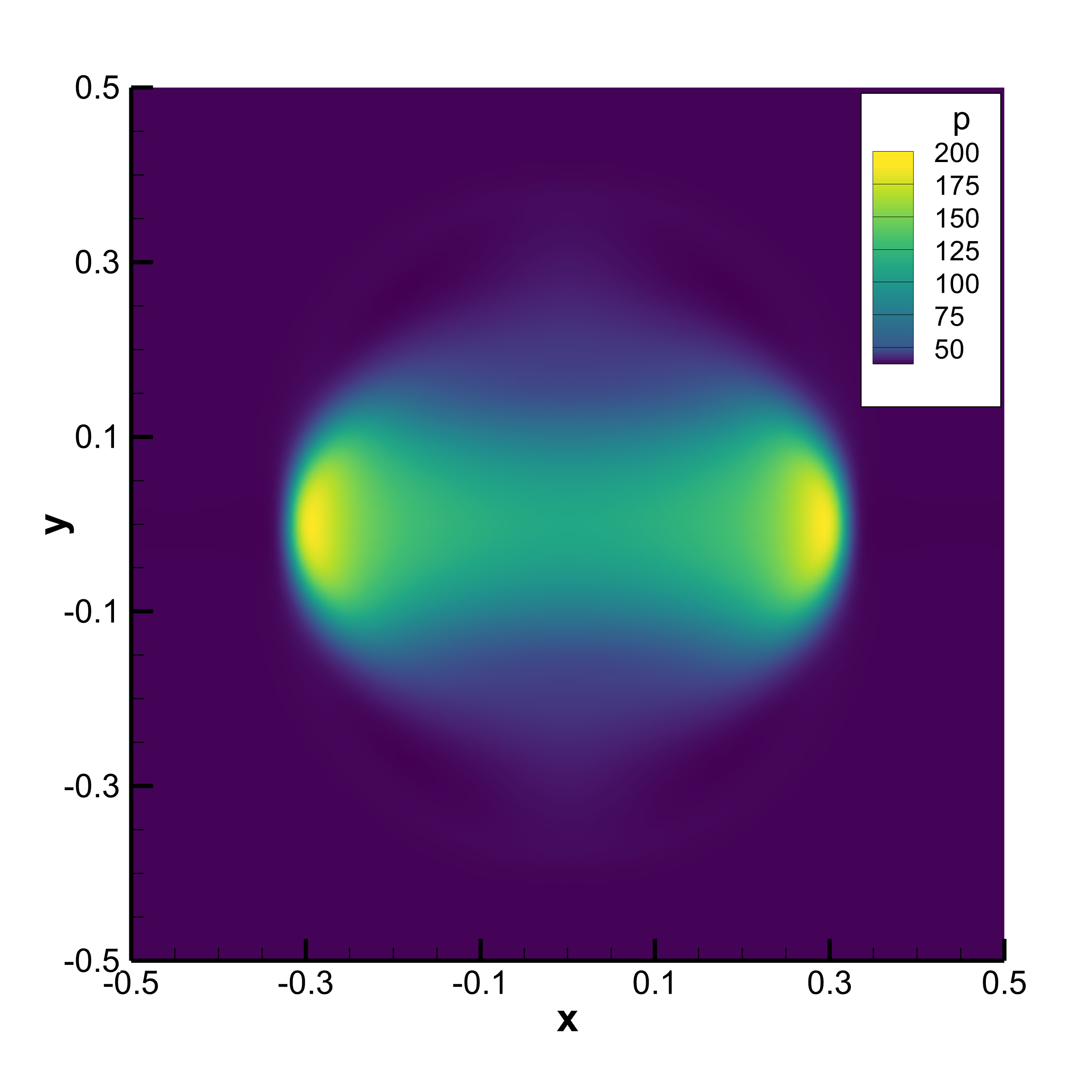}\\[-1mm]
\footnotesize{(b) Pressure}
\end{minipage}\hfill
\begin{minipage}[t]{0.32\textwidth}\centering
\includegraphics[width=\linewidth]{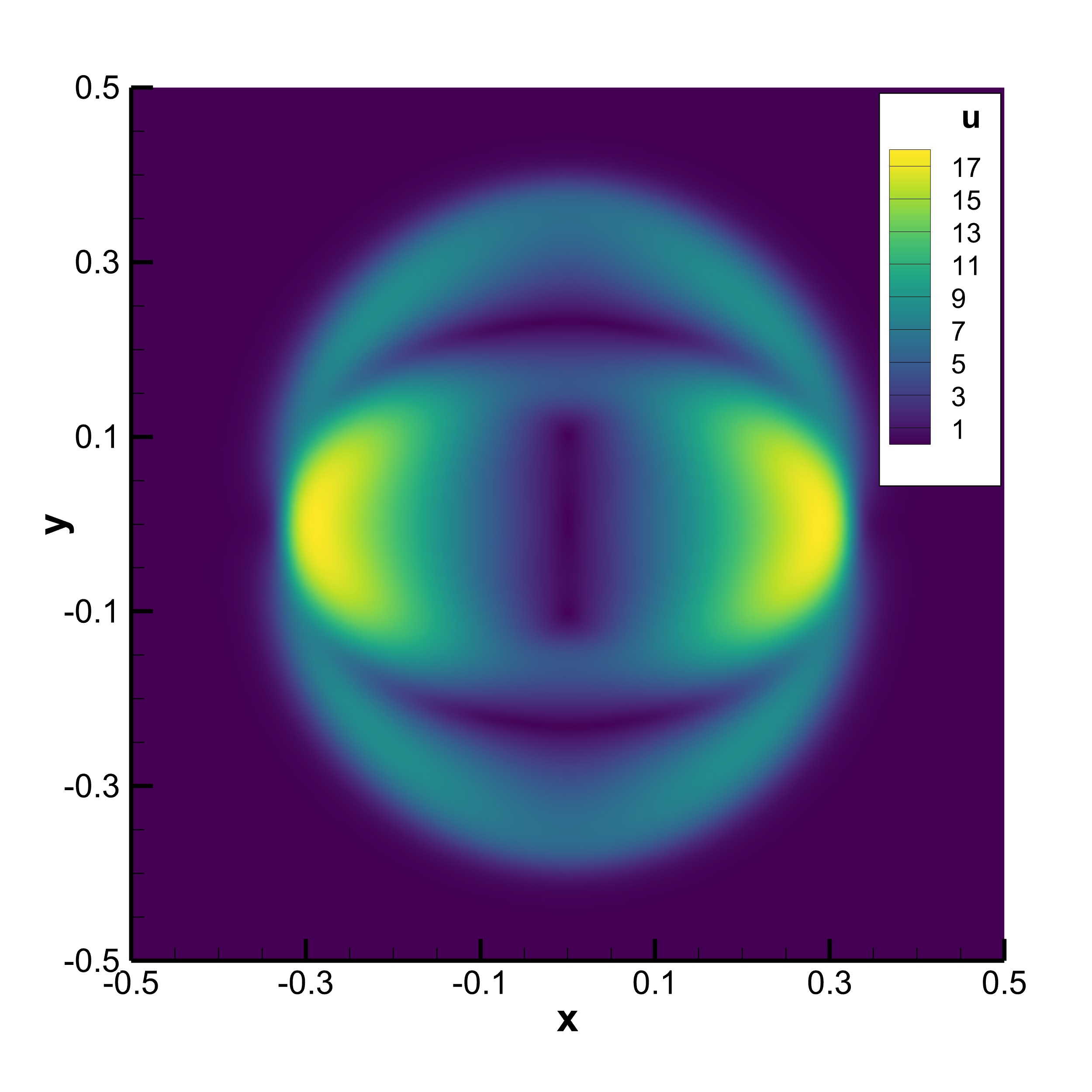}\\[-1mm]
\footnotesize{(c) Total velocity}
\end{minipage}

\medskip

\begin{minipage}[t]{0.32\textwidth}\centering
\includegraphics[width=\linewidth]{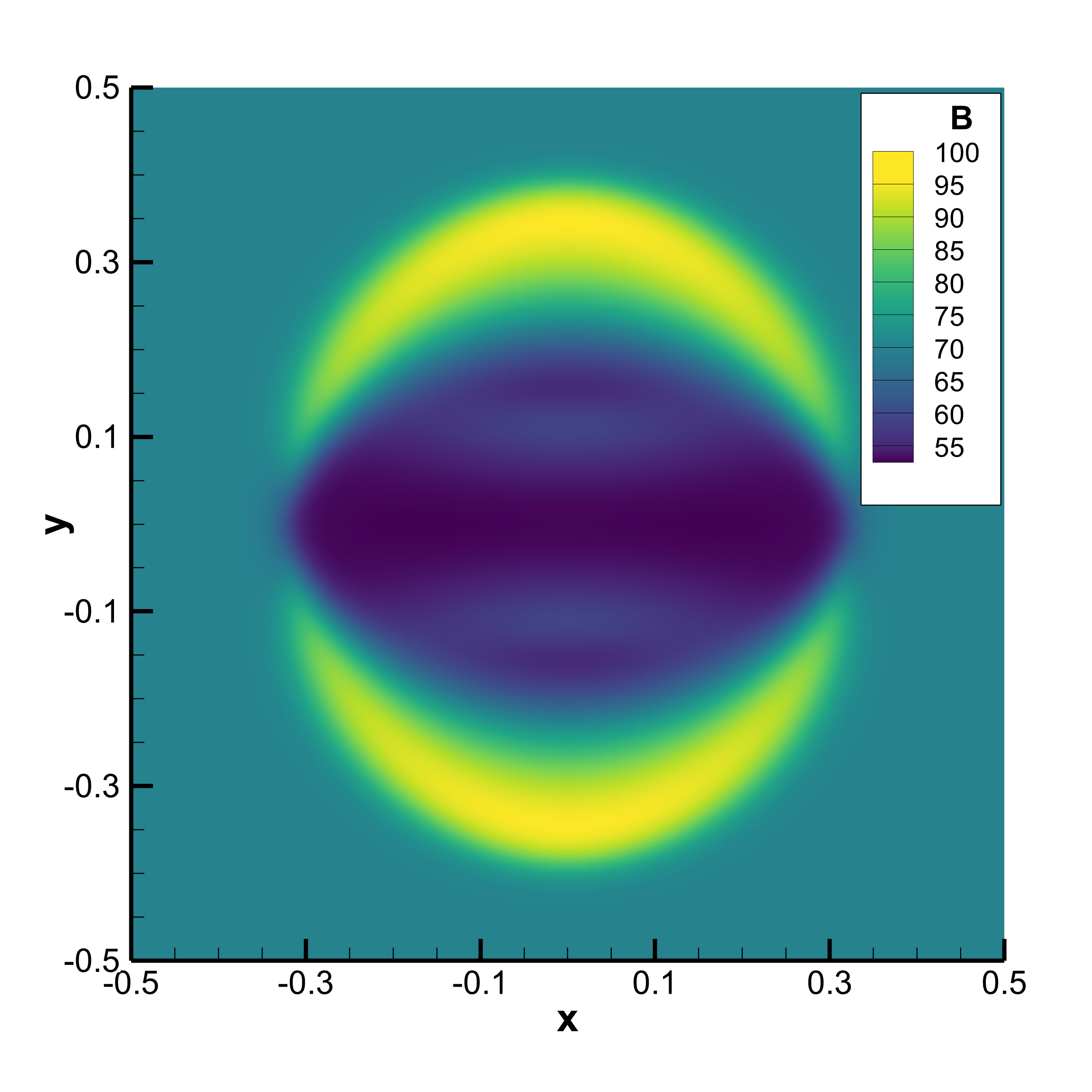}\\[-1mm]
\footnotesize{(d) Total magnetic field}
\end{minipage}
\begin{minipage}[t]{0.32\textwidth}\centering
\includegraphics[width=\linewidth]{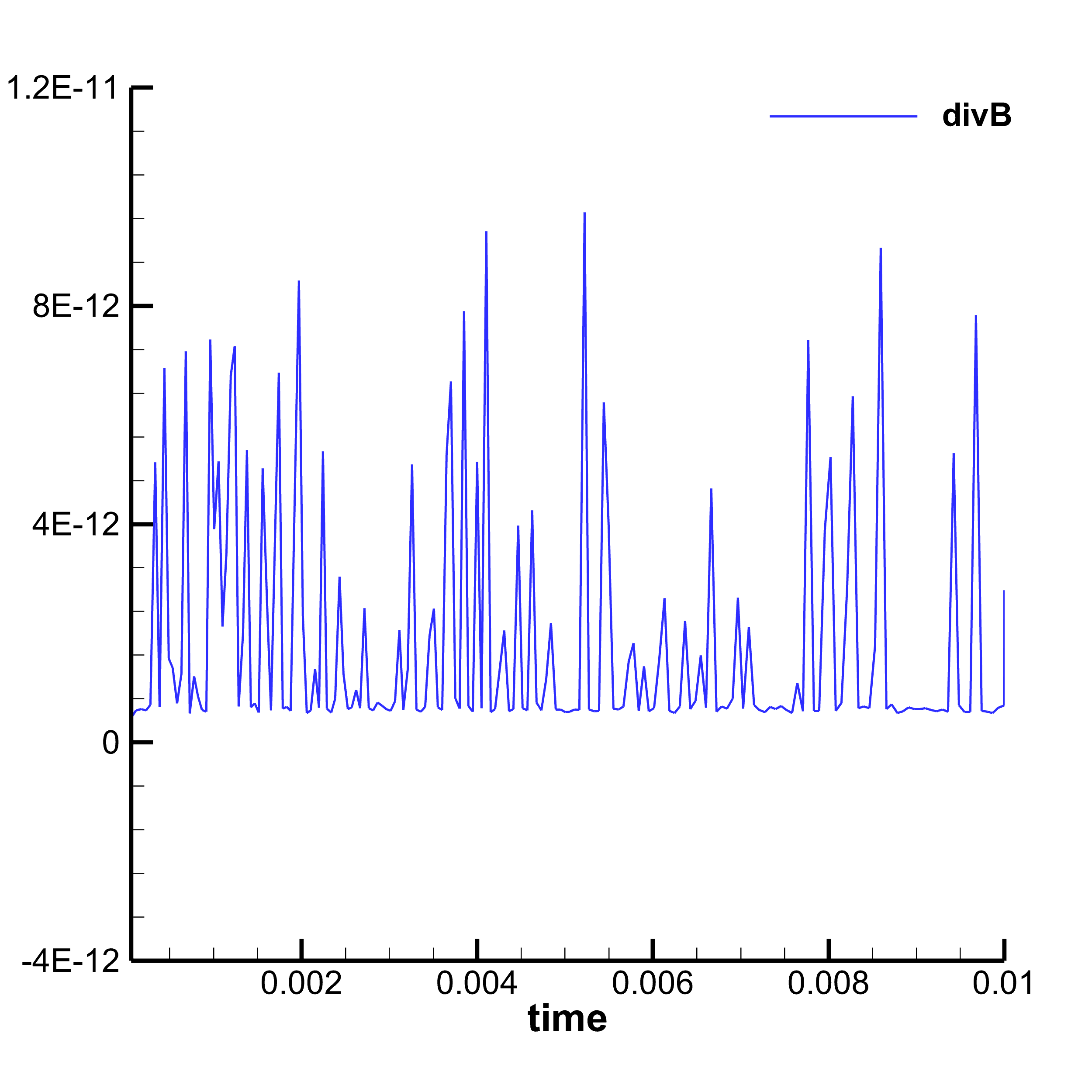}\\[-1mm]
\footnotesize{(e) Time evolution of $\nabla \cdot \mathbf{B}$}
\end{minipage}
\caption{Numerical solution of the \textit{MHD Blast wave} at the final time $t_f = 0.01$ obtained using our involution preserving scheme on a $N_x = N_y = 200$ mesh.}
\label{fig:mhd_blast}
\end{figure}

The MHD blast wave problem introduced in \cite{balsara1999staggered} represents a challenging test since it contains a strong magnetic field and a jump in the pressure of four orders of magnitude. The initial density is constant $\rho = 1$, the initial velocity field is set to zero $(u, v) = (0, 0)$ and the magnetic field is $\bm{B} = (70, 0)$. The computational domain is $\Omega = [-0.5, 0.5]^2$. Inside the circular region $r < 0.1$, with $r = \sqrt{x^2 + y^2}$, the pressure is set to $p = 1000$ whereas outside is assigned with $p = 0.1$. Thus initially, the Alfvén waves are faster than the sound speed inside the circular domain. Here we set $\gamma = 1.4$ and the final time is chosen to be $t_f = 0.01$. The numerical solution at the final time is reported in Figure~\ref{fig:mhd_blast}. The maximum discrete divergence over time is $9.71\times10^{-12}$ and its final value is $2.79\times10^{-12}$.

\begin{figure}[t]
\centering
\begin{minipage}[t]{0.32\textwidth}\centering
\includegraphics[width=\linewidth]{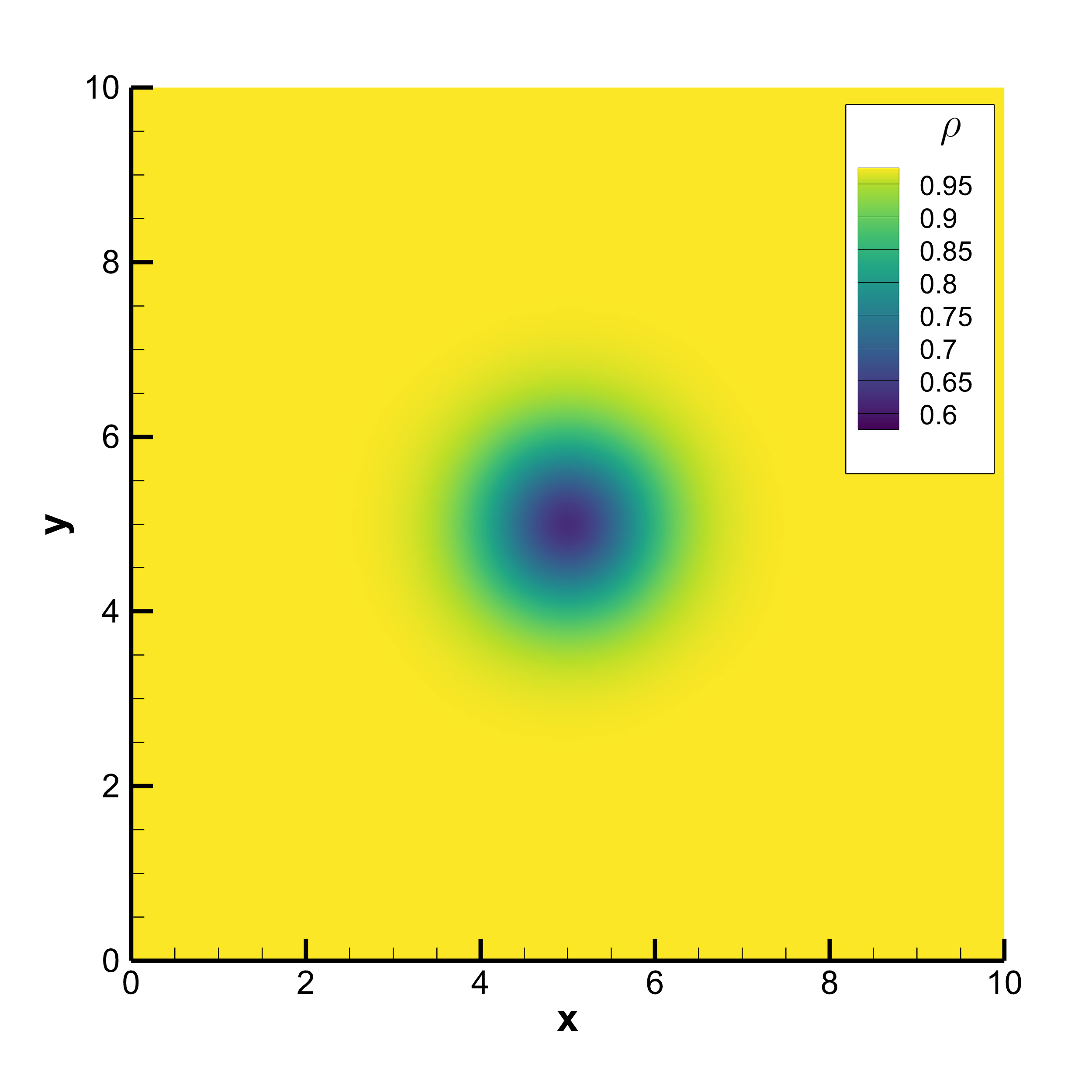}\\[-1mm]
\footnotesize{(a) Density}
\end{minipage}\hfill
\begin{minipage}[t]{0.32\textwidth}\centering
\includegraphics[width=\linewidth]{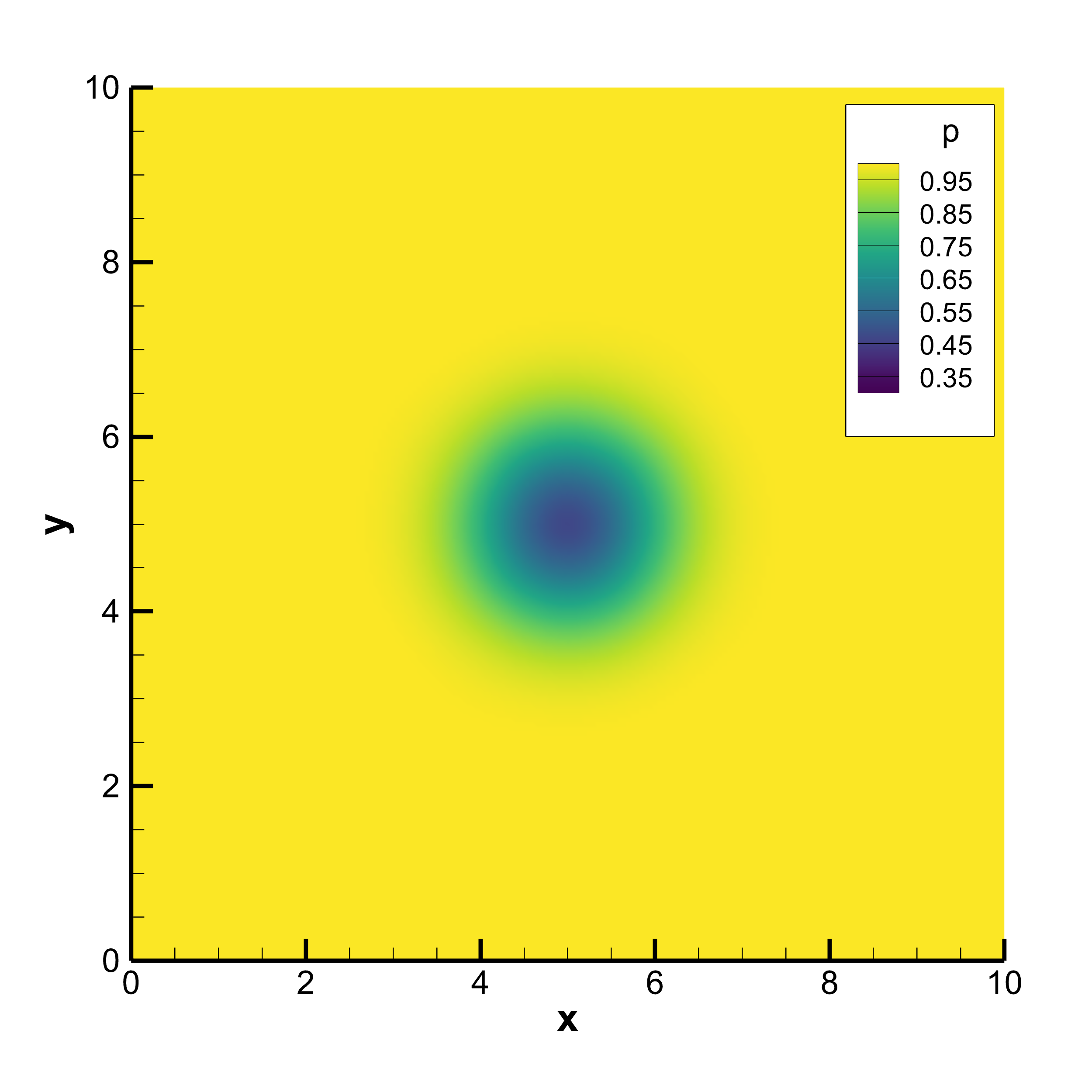}\\[-1mm]
\footnotesize{(b) Pressure}
\end{minipage}\hfill
\begin{minipage}[t]{0.32\textwidth}\centering
\includegraphics[width=\linewidth]{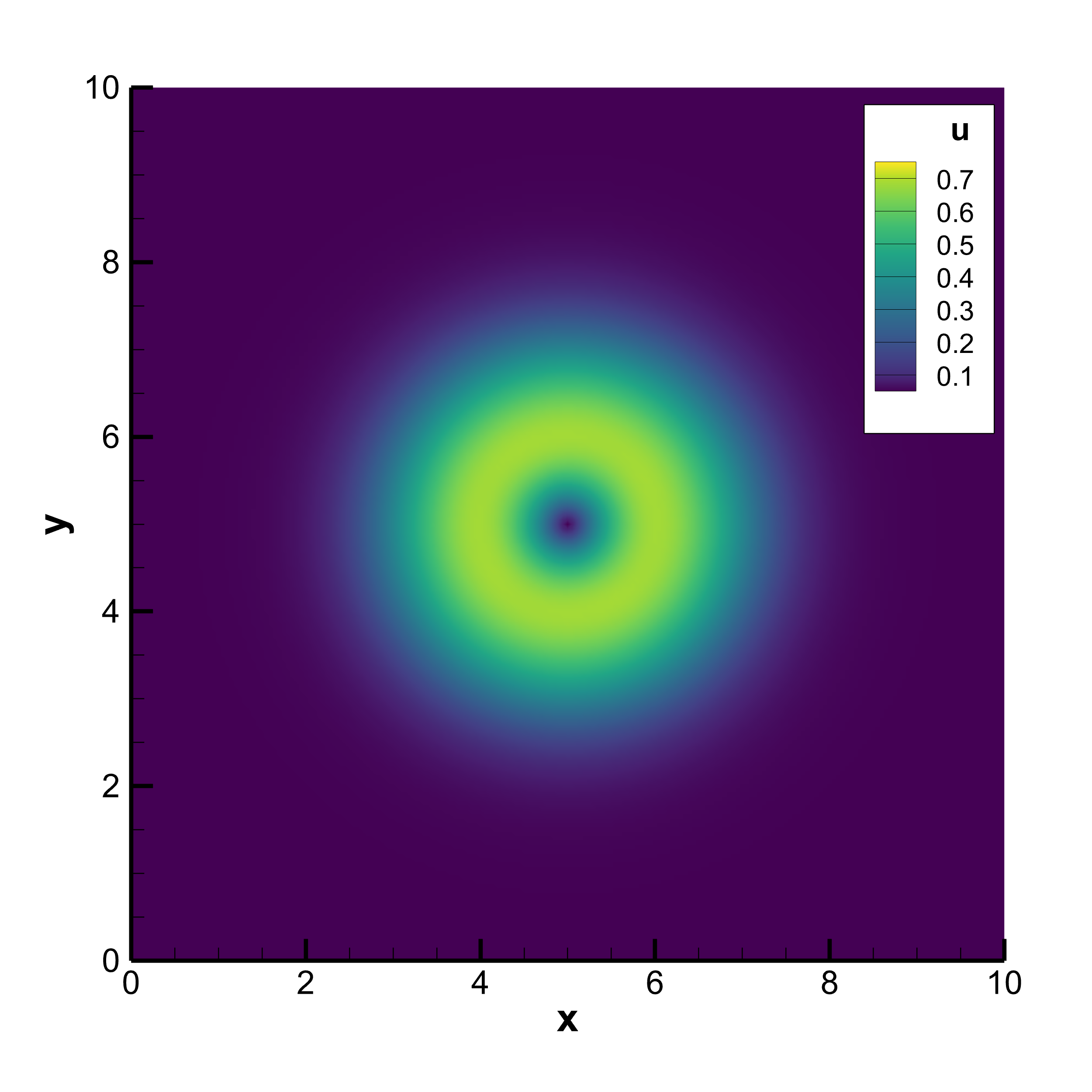}\\[-1mm]
\footnotesize{(a) Total velocity}
\end{minipage}

\medskip

\begin{minipage}[t]{0.32\textwidth}\centering
\includegraphics[width=\linewidth]{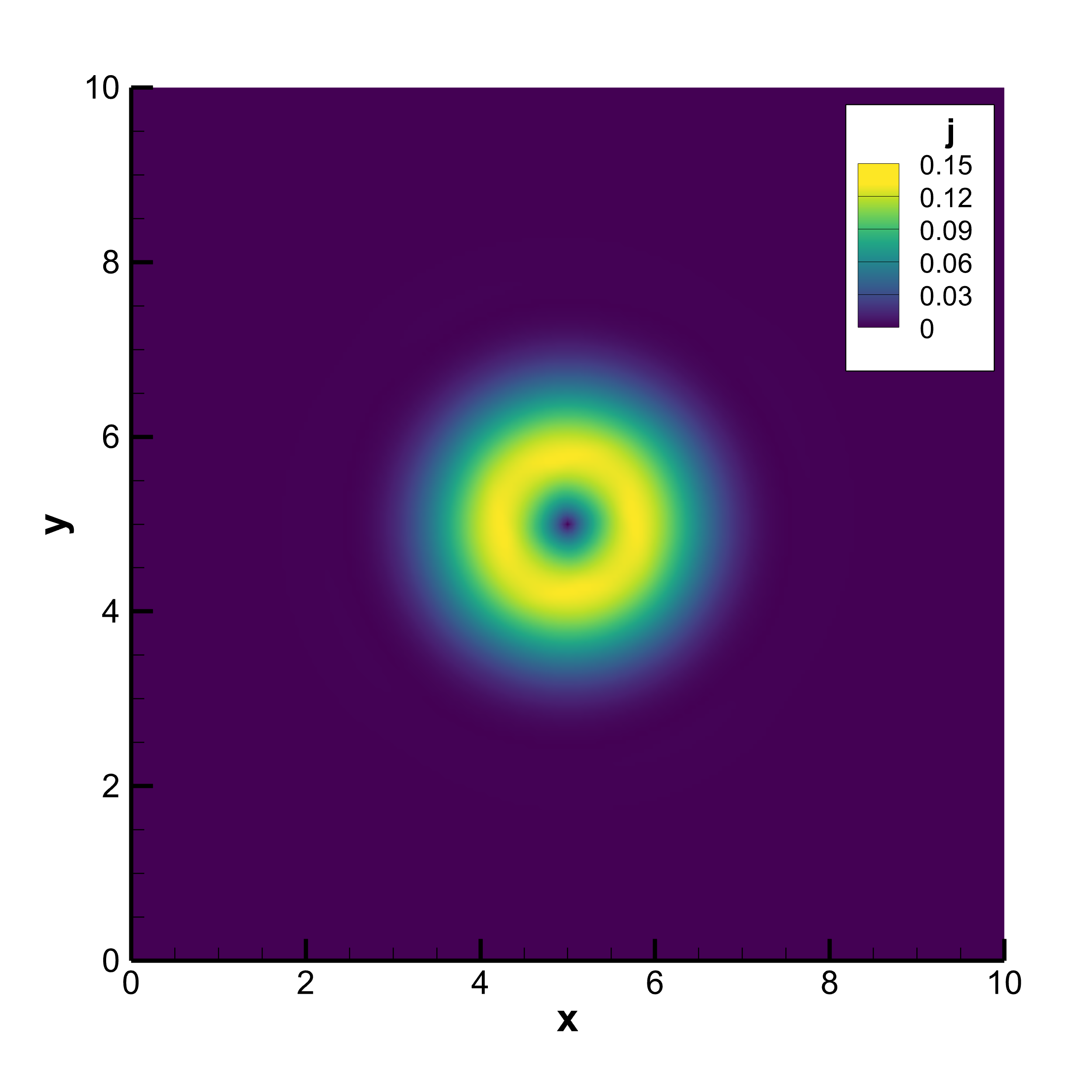}\\[-1mm]
\footnotesize{(d) Total thermal impulse}
\end{minipage}
\begin{minipage}[t]{0.32\textwidth}\centering
\includegraphics[width=\linewidth]{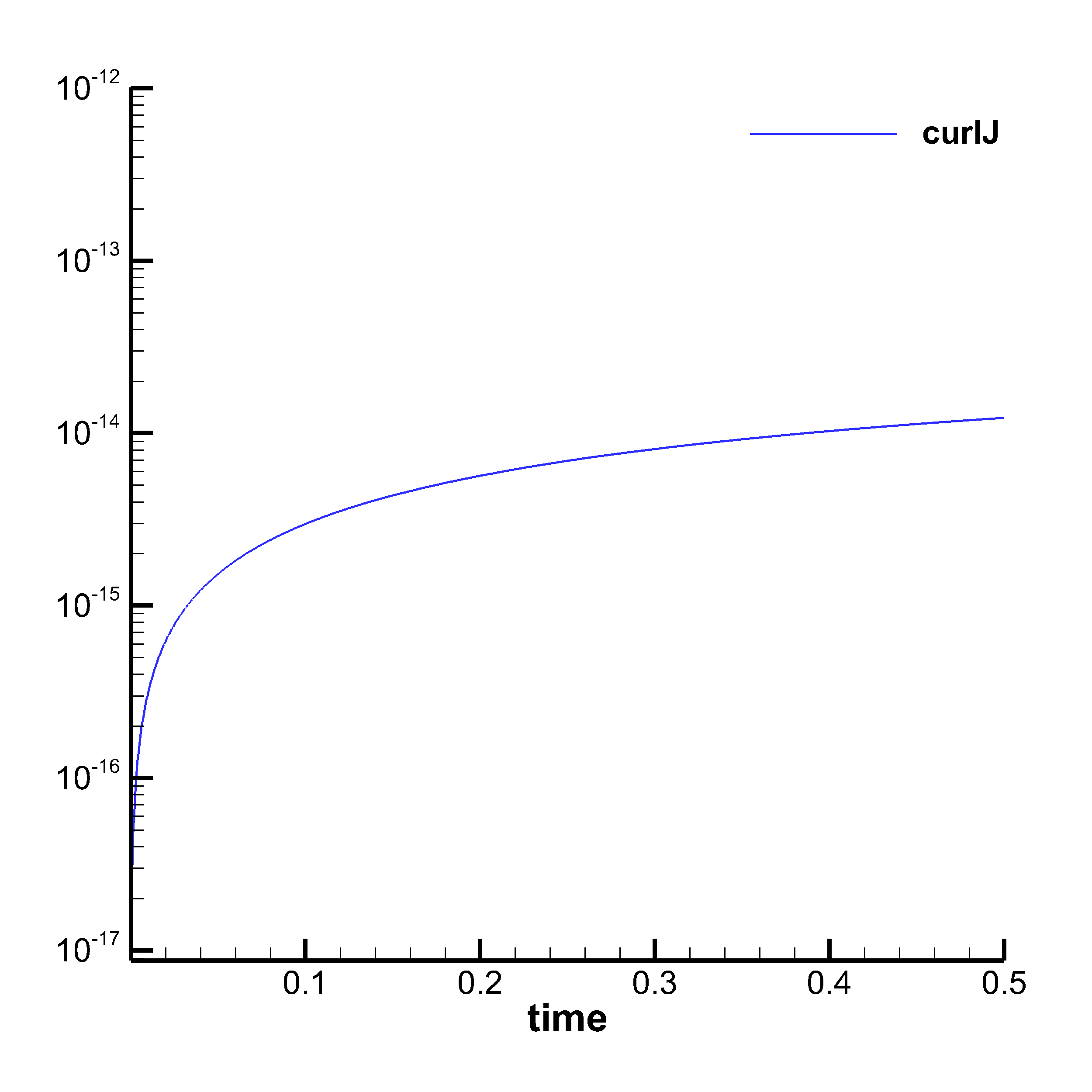}\\[-1mm]
\footnotesize{(e) Time evolution of $\nabla \times \mathbf{j}$}
\end{minipage}
\caption{Numerical solution of the \textit{Euler-Heat vortex with dissipation} at the final time $t_f = 0.5$ obtained using our involution preserving scheme on a $N_x = N_y = 200$ mesh.}
\label{fig:eh_dissipative_vortex}
\end{figure}

\subsubsection{Euler-Heat numerical convergence studies}

To assess the accuracy, we consider a stationary manufactured solution for the Euler--Heat model. The computational domain is the square $\Omega = [0, 10]^2$. We use the centred coordinates $\widetilde{x}:=x-5$, $\widetilde{y}:=y-5$, and $r:=\sqrt{\widetilde{x}^2+\widetilde{y}^2}$. The fluid is assigned a unit density and pressure with the following velocity and scalar potential $\varphi$:
\begin{equation}
	\begin{aligned}
		\begin{pmatrix}
			u \\
			v
		\end{pmatrix}
		&=
		\begin{pmatrix}
			-\dfrac{\beta e^{\frac{1}{2}(1-r^2)}}{2\pi}\,\widetilde{y} \\
			\dfrac{\beta e^{\frac{1}{2}(1-r^2)}}{2\pi}\,\widetilde{x}
		\end{pmatrix},
		\qquad
		\varphi =
		-\dfrac{\lambda e^{\frac{1}{2}(1-r^2)}}{2\pi},
		\qquad \beta,\lambda\in\mathbb{R}.
	\end{aligned}
\end{equation}
The thermal impulse is initialized by $\bm{j}=\nabla\varphi$, which is curl-free in the continuous sense. We choose $\beta=1$ and $\lambda=1$. The stationary exact state $\bm{q}^{\mathrm{ex}}$ is maintained by the manufactured source vector
\begin{align}
  \bm{S}^{\mathrm{ex}}(\bm{x})
  :=\nabla\cdot\mathds{F}^{\scriptscriptstyle\mathrm{EH}}
  (\bm{q}^{\mathrm{ex}}(\bm{x})),
  \label{eq:eh_manufactured_source}
\end{align}
which is evaluated consistently for all six equations of the model.

The simulation is run on a sequence of successively refined meshes until the final time $t_f = 0.1$. Both the $L^1$ and the $L^2$ errors are reported in Table \ref{tab:eh_convergence}, showing first-order accuracy, as expected.

\begin{table}[htbp]
	\centering
	\caption{Euler--Heat convergence on moving meshes at $t_f=0.1$.
		$L^1$ and $L^2$ errors and experimental orders of convergence for the primitive variables.}
	\label{tab:eh_convergence}
	\scriptsize
	\setlength{\tabcolsep}{3pt}
	
	\begin{tabular}{c cc cc cc cc cc cc}
		\hline
		$N_x=N_y$
		& \multicolumn{2}{c}{$\rho$}
		& \multicolumn{2}{c}{$u$}
		& \multicolumn{2}{c}{$v$}
		& \multicolumn{2}{c}{$p$}
		& \multicolumn{2}{c}{$j_x$}
		& \multicolumn{2}{c}{$j_y$} \\
		\hline
		
		\multicolumn{13}{c}{$L^1$ errors} \\
		\hline
		
		20
		& $3.126\mathrm{e}{-3}$ & -
		& $5.034\mathrm{e}{-3}$ & -
		& $5.034\mathrm{e}{-3}$ & -
		& $3.113\mathrm{e}{-3}$ & -
		& $5.027\mathrm{e}{-3}$ & -
		& $5.027\mathrm{e}{-3}$ & - \\
		
		30
		& $2.126\mathrm{e}{-3}$ & 0.957
		& $3.447\mathrm{e}{-3}$ & 0.940
		& $3.447\mathrm{e}{-3}$ & 0.940
		& $2.102\mathrm{e}{-3}$ & 0.974
		& $3.436\mathrm{e}{-3}$ & 0.945
		& $3.436\mathrm{e}{-3}$ & 0.945 \\
		
		40
		& $1.657\mathrm{e}{-3}$ & 0.875
		& $2.668\mathrm{e}{-3}$ & 0.901
		& $2.668\mathrm{e}{-3}$ & 0.901
		& $1.633\mathrm{e}{-3}$ & 0.887
		& $2.631\mathrm{e}{-3}$ & 0.938
		& $2.631\mathrm{e}{-3}$ & 0.938 \\
		
		50
		& $1.352\mathrm{e}{-3}$ & 0.926
		& $2.186\mathrm{e}{-3}$ & 0.905
		& $2.186\mathrm{e}{-3}$ & 0.905
		& $1.331\mathrm{e}{-3}$ & 0.930
		& $2.145\mathrm{e}{-3}$ & 0.930
		& $2.145\mathrm{e}{-3}$ & 0.930 \\
		
		60
		& $1.153\mathrm{e}{-3}$ & 0.888
		& $1.859\mathrm{e}{-3}$ & 0.904
		& $1.859\mathrm{e}{-3}$ & 0.904
		& $1.128\mathrm{e}{-3}$ & 0.921
		& $1.810\mathrm{e}{-3}$ & 0.943
		& $1.810\mathrm{e}{-3}$ & 0.943 \\
		
		\hline
		\multicolumn{13}{c}{$L^2$ errors} \\
		\hline
		
		20
		& $7.628\mathrm{e}{-3}$ & -
		& $1.190\mathrm{e}{-2}$ & -
		& $1.190\mathrm{e}{-2}$ & -
		& $7.356\mathrm{e}{-3}$ & -
		& $1.181\mathrm{e}{-2}$ & -
		& $1.181\mathrm{e}{-2}$ & - \\
		
		30
		& $5.302\mathrm{e}{-3}$ & 0.903
		& $8.353\mathrm{e}{-3}$ & 0.879
		& $8.353\mathrm{e}{-3}$ & 0.879
		& $5.025\mathrm{e}{-3}$ & 0.946
		& $8.259\mathrm{e}{-3}$ & 0.888
		& $8.259\mathrm{e}{-3}$ & 0.888 \\
		
		40
		& $4.197\mathrm{e}{-3}$ & 0.821
		& $6.492\mathrm{e}{-3}$ & 0.886
		& $6.492\mathrm{e}{-3}$ & 0.886
		& $3.949\mathrm{e}{-3}$ & 0.847
		& $6.369\mathrm{e}{-3}$ & 0.914
		& $6.369\mathrm{e}{-3}$ & 0.914 \\
		
		50
		& $3.454\mathrm{e}{-3}$ & 0.887
		& $5.375\mathrm{e}{-3}$ & 0.858
		& $5.375\mathrm{e}{-3}$ & 0.858
		& $3.241\mathrm{e}{-3}$ & 0.899
		& $5.240\mathrm{e}{-3}$ & 0.887
		& $5.240\mathrm{e}{-3}$ & 0.887 \\
		
		60
		& $2.895\mathrm{e}{-3}$ & 0.982
		& $4.590\mathrm{e}{-3}$ & 0.880
		& $4.590\mathrm{e}{-3}$ & 0.880
		& $2.715\mathrm{e}{-3}$ & 0.985
		& $4.450\mathrm{e}{-3}$ & 0.910
		& $4.450\mathrm{e}{-3}$ & 0.910 \\
		
		\hline
	\end{tabular}
\end{table}


\subsubsection{Euler-Heat vortex problem with dissipation}

The computational domain is $\Omega = [0, 10]^2$. The initial condition reads
\begin{equation}
	\begin{aligned}
		\rho = (1 + \delta\theta)^{\frac{1}{\gamma - 1}},
		\qquad
		\begin{pmatrix}
			u\\
			v
		\end{pmatrix}
		= 
			\frac{5}{2\pi} e^{\frac{1 - r^2}{2}}
			\begin{pmatrix}
				5 - y\\
				x - 5
			\end{pmatrix},
		\qquad
		\delta\theta = - \frac{25(\gamma - 1)}{8\gamma \pi^2}e^{1-r^2},
	\end{aligned}
\end{equation}
with $p = \rho ^ \gamma$ and $r = \sqrt{(x-5)^2 + (y-5)^2}$. The thermal impulse is initially $\bm{j}=\bm{0}$. We set the hyperbolic heat parameter to $\kappa = 10^{-2}$ and the target Fourier conductivity to $K = 10^{-3}$. The constant relaxation time is then $\tau = K/\kappa^2$, hence preserving a zero curl under the relaxation source $-\bm{j}/\tau$ in the evolution equations for $\bm{j}$. The simulation is run until the final time $t_f = 0.5$. The resulting vortex and the discrete curl are shown in Figure~\ref{fig:eh_dissipative_vortex}.

This test is advanced by a symmetric source splitting: one half of the exact relaxation flow, the homogeneous constrained ALE step, and a second half relaxation flow. For the two components of $\bm{j}$, the centered nodal flux is used, whereas the remaining equations retain the stabilized numerical flux. On the initially Cartesian mesh, both this compatible flux and the scalar relaxation map preserve the discrete curl. Consequently, the constrained mesh problem is still posed at every Runge--Kutta stage, but its admissible minimum is the zero displacement and the computed grid velocity remains zero up to roundoff. This behavior is noteworthy: the method does not introduce mesh motion when the spatial discretization already satisfies the involution. The experiment is therefore complementary to the preceding convergence test, in which a nonzero ALE motion is required to enforce the curl constraint. The maximum and final raw discrete curl are both $1.23\times10^{-14}$. 

\subsection{Geometric Conservation Law}

The proposed ALE discretization is also assessed in terms of the discrete Geometric
Conservation Law (GCL), which guarantees that the numerical evolution of the cell
volumes is exactly consistent with the prescribed mesh motion. This property is
essential for the preservation of uniform solutions on moving meshes.

For each control volume $\omega_c$, the continuous GCL reads
\begin{equation}
	\frac{\mathrm d}{\mathrm dt} |\omega_c|
	=
	\sum_{p\in\mathcal P(c)}
	\bm{\mathcal w}_p\cdot\bm{\ell}_{pc},
\end{equation}
which is the geometric counterpart of the ALE finite-volume formulation.

To verify the discrete GCL, we compare the cell volume obtained after updating the
mesh coordinates through the Runge--Kutta scheme with the volume predicted along
the same piecewise-linear displacement path. The corner vectors are evaluated at
the midpoint of each segment, as in \eqref{eq:segment_midpoint_gcl_flux}, rather
than simply on the physical Runge--Kutta stage geometry. Denoting by
$|\omega_c|_{\mathrm{geom}}^{n+1}$ the cell volume computed from the updated mesh,
the GCL residual is defined as
\begin{equation}
	R_{\mathrm{GCL},c}
	=
	|\omega_c|_{\mathrm{geom}}^{n+1}
	-
	|\omega_c|^{n}
	-
	\Delta t
	\sum_{s=1}^{N_{\mathrm{RK}}}
	\mathsf b_s
	G_{c,s}^{\mathrm{mid}},
\end{equation}
where
\begin{equation}
	G_{c,s}^{\mathrm{mid}}
	=
	\sum_{p\in\mathcal P(c)}
	\bm{\mathcal w}_{p}^{(s)}
	\cdot
	\bm{\ell}_{pc}(\bm{X}_{s-1/2})
\end{equation}
is the exact swept-volume flux for the $s$-th displacement segment. This
distinction is essential because polygonal area is quadratic in the coordinates.

The reported diagnostic is the normalized residual
\begin{equation}
	\varepsilon_{\mathrm{GCL}}
	=
	\max_c
	\frac{|R_{\mathrm{GCL},c}|}
	{\max\!\left(
		|\omega_c^{n}|,
		|\omega_c^{n+1}|,
		\delta_{\mathrm{m}}
		\right)},
\end{equation}
with $\delta_{\mathrm{m}}$ denoting the machine precision of the floating-point arithmetic.

\begin{figure}[t]
	\centering
	\begin{minipage}[t]{0.32\textwidth}\centering
		\includegraphics[width=\linewidth]{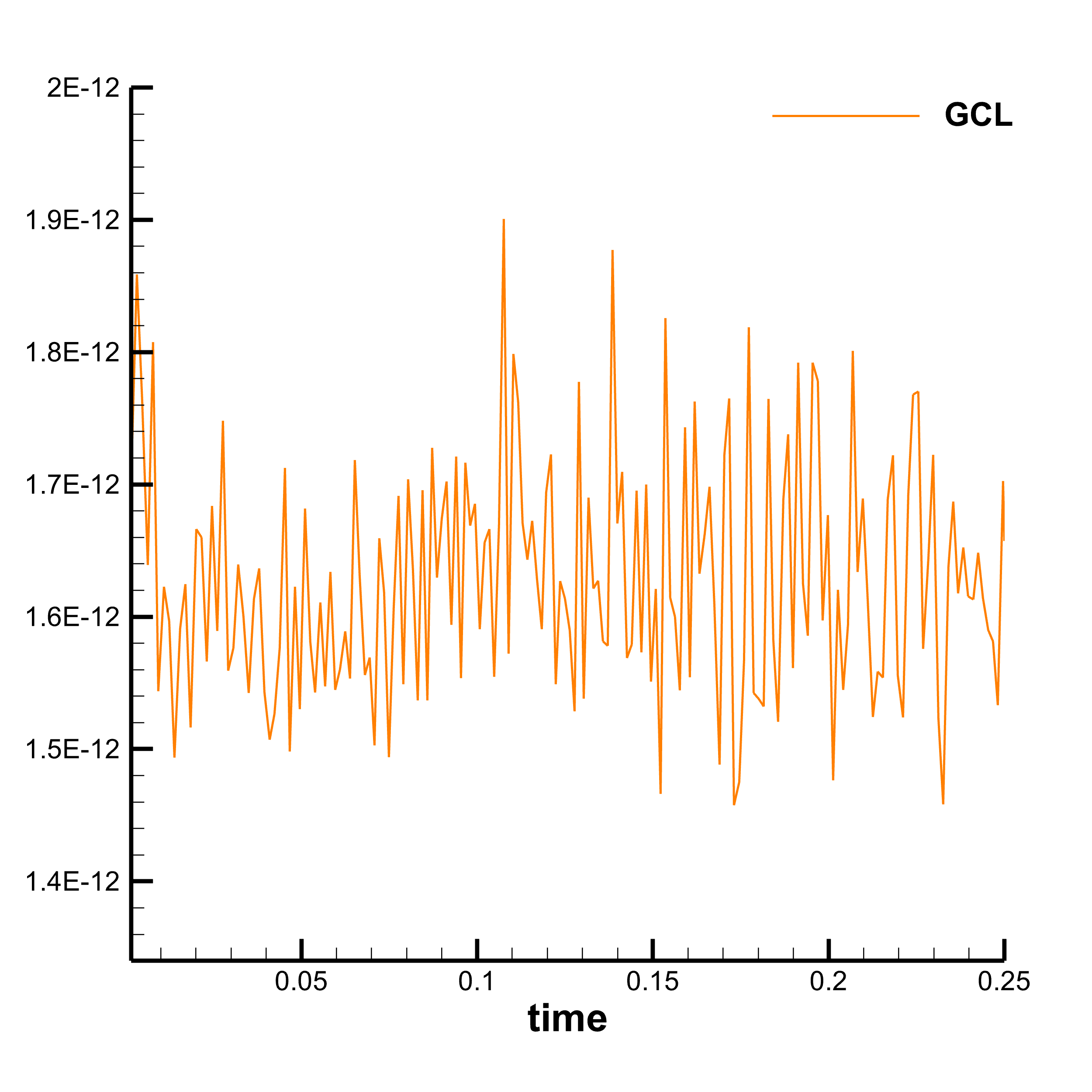}\\[-1mm]
		\footnotesize{(a) \textit{Ideal-MHD Rotor}}
	\end{minipage}\hfill
	\begin{minipage}[t]{0.32\textwidth}\centering
		\includegraphics[width=\linewidth]{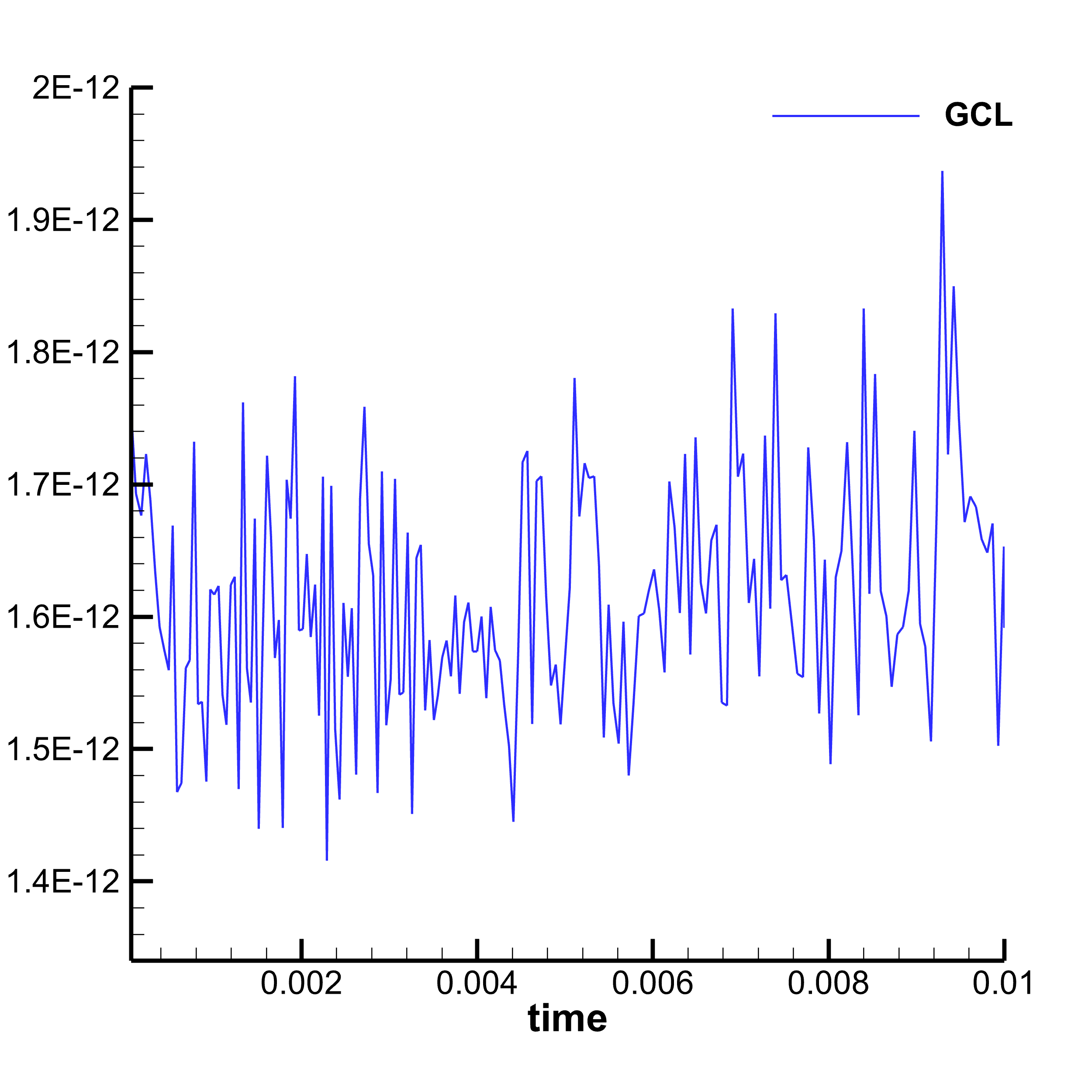}\\[-1mm]
		\footnotesize{(b) \textit{Ideal-MHD blast wave}}
	\end{minipage}\hfill
	\begin{minipage}[t]{0.32\textwidth}\centering
		\includegraphics[width=\linewidth]{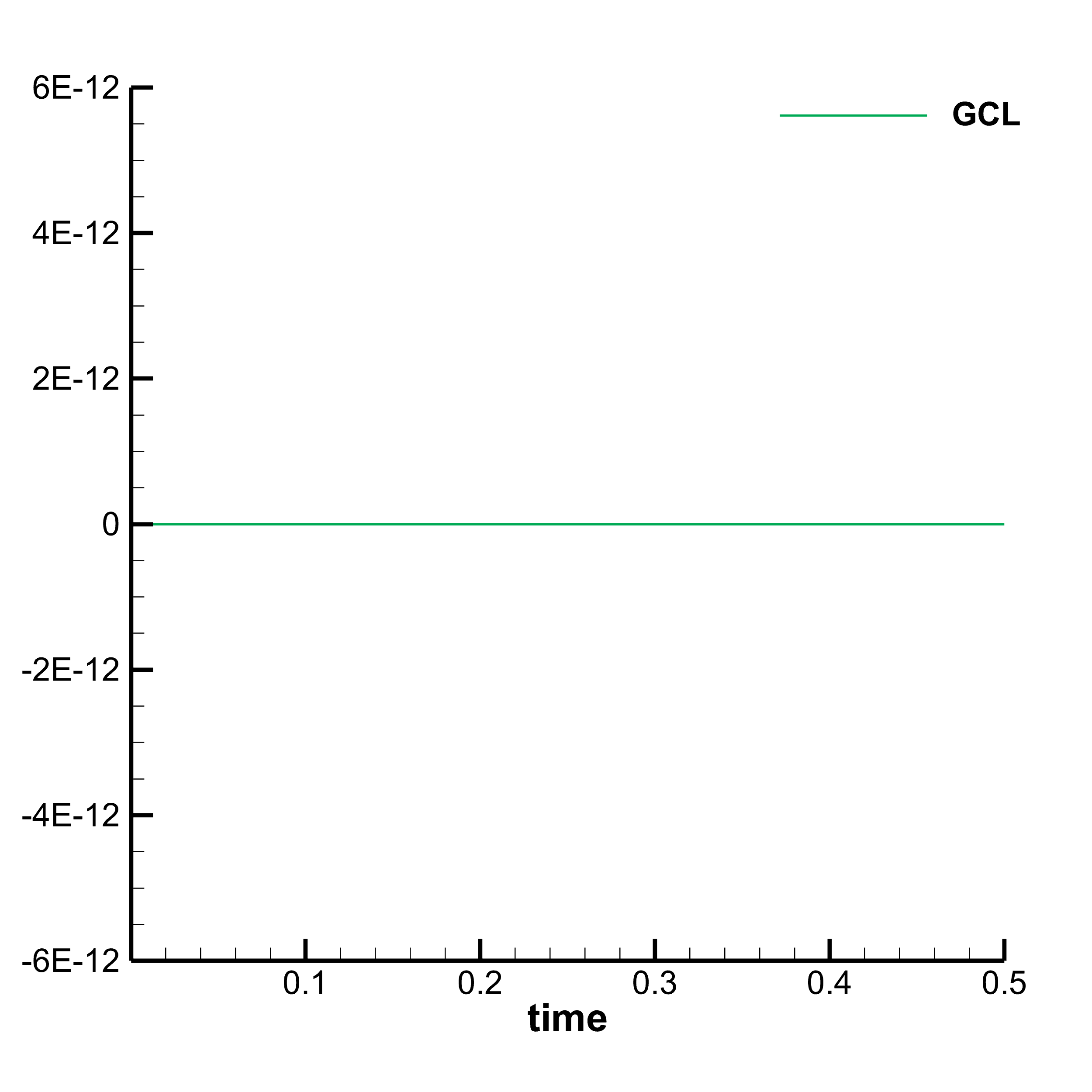}\\[-1mm]{}
		\footnotesize{(b) \textit{Dissipative Euler-Heat vortex}}
	\end{minipage}
	\caption{Time evolution of the maximum normalized GCL residual 
		$\varepsilon_{\mathrm{GCL}}$ for the different test cases.}
	\label{GCL plots}
\end{figure}

Furthermore, the same verification is performed at every Runge--Kutta stage,
ensuring that the geometry evolution is consistent throughout the multistage
integration. In all numerical experiments, both the stage-wise and the final GCL
residuals remain at the level of machine precision, confirming the exact
satisfaction of the discrete GCL, as shown in Figure \ref{GCL plots}.
As discussed previously, for the Euler--Heat dissipation vortex test, the constrained mesh problem admits the zero displacement as its solution and the resulting grid velocity remains zero up to roundoff. Consequently, the mesh remains stationary and both the stage-wise and final GCL residuals are identically zero.

\section{Conclusions} \label{sec.concl}
In this work, we have introduced an involution-preserving ALE finite volume framework for the numerical approximation of first-order hyperbolic systems on moving meshes. The main idea consists in exploiting the additional freedom provided by the ALE formulation to construct the mesh motion consistently with the differential involutions associated to the governing equations.

The method has been formulated in a cell-centered setting on initially quadrilateral meshes. Corner-based geometrical quantities are devised to define the discrete differential operators, while all geometrical contributions are consistently updated along the mesh evolution. The resulting constrained mesh motion has been incorporated into a second-order Runge--Kutta time discretization, ensuring consistency between the intermediate stages, the evolving geometry, and the imposed differential constraint. The formulation also satisfies the Geometric Conservation Law, preventing the introduction of spurious errors emerging from the mesh motion.

The proposed approach has been investigated for two different hyperbolic systems. For the ideal MHD equations, the mesh motion is constructed in connection with the divergence-free involution of the magnetic field, whereas for the hyperbolic Euler--Heat system the corresponding construction is associated with the curl-free involution of the thermal impulse vector. The numerical results confirm the expected accuracy of the scheme and show that the relevant differential constraints remain satisfied throughout the evolution on moving meshes at the fully discrete level up to machine precision. Classical benchmark problems further demonstrate the capability of the method to handle non-trivial flow configurations while maintaining the prescribed structural properties.

Future developments will concern the extension of the proposed framework to higher-order spatial and temporal discretizations and to polygonal moving meshes. The satisfaction of multiple involutions at the same time will also be the subject of investigations, as well as the application of our novel strategy to other hyperbolic systems like hyperelasticity.

\section*{Acknowledgments}
The authors would like to thank the framework offered by the SUN HYPE 2026 event organized in Crete in June 2026, where the work presented in this paper was conceived and its initial development took place. In particular, we thank the Technical University of Crete for hosting the event and the financial support received from the European Union through the ERC Starting Grant ALcHyMiA (grant agreement No. 101114995). Views and opinions expressed are however those of the author(s) only and do not necessarily reflect those of the European Union or the European Research Council Executive Agency. Neither the European Union nor the granting authority can be held responsible for them.

\newpage

\bibliographystyle{siamplain}
\bibliography{biblio}

@Article{Maire2009,
  author  = {{P.-H.} Maire},
  title   = {A high-order cell-centered Lagrangian scheme for two-dimensional compressible fluid flows on unstructured meshes.},
  journal = {J. Comput. Phys.},
  year    = {2009},
  volume  = {228},
  pages   = {2391--2425},
}

@article{KnuppRezoning,
  author = {{P.M.} Knupp},
  title = {Achieving finite element mesh quality via optimization of the Jacobian matrix norm and associated quantities. Part II -- A framework for volume mesh optimization and the condition number of the Jacobian matrix.},
  journal = {Int. J. Numer. Meth. Engng.},
  year = {2000},
  volume = {48},
  pages = {1165 - 1185}
}

@article{gottlieb2001strong,
  title   = {Strong stability-preserving high-order
             time discretization methods},
  author  = {Gottlieb, Sigal and Shu, Chi-Wang
             and Tadmor, Eitan},
  journal = {SIAM Review},
  volume  = {43},
  number  = {1},
  pages   = {89--112},
  year    = {2001}
}

@article{BoscheriGPRGCL,
author = {W. Boscheri and R. Loub\`ere and J.P. Braeunig and P.-H. Maire},
title = {A geometrically and thermodynamically compatible finite volume scheme for continuum mechanics on unstructured polygonal meshes},
journal = {J. Comput. Phys.},
volume = {507},
pages = {112957},
year = {2024},
}

@article{DhaouadiGavrilyuk2024,
  author  = {F. Dhaouadi and S. Gavrilyuk},
  title   = {An {Eulerian} hyperbolic model for heat transfer derived via
             {Hamilton}'s principle: analytical and numerical study},
  journal = {Proc. R. Soc. A},
  volume  = {480},
  number  = {2283},
  pages   = {20230440},
  year    = {2024},
}

@article{BoscheriDhaouadi2026,
  author  = {W. Boscheri and F. Dhaouadi},
  title   = {Structure preserving finite volume schemes on {Voronoi} grids:
             curl involution, asymptotic limit and thermodynamics},
  journal = {J. Sci. Comput.},
  volume  = {107},
  pages   = {80},
  year    = {2026},
  doi     = {10.1007/s10915-026-03251-6},
}

@article{balsara1999staggered,
  title={A staggered mesh algorithm using high order Godunov fluxes to ensure solenoidal magnetic fields in magnetohydrodynamic simulations},
  author={Balsara, Dinshaw S and Spicer, Daniel S},
  journal={Journal of Computational Physics},
  volume={149},
  number={2},
  pages={270--292},
  year={1999},
  publisher={Elsevier}
}

@article{Tonon2021,
  author  = {P. Tonon and R.A.K. Sanches and K. Takizawa and T.E. Tezduyar},
  title   = {A linear-elasticity-based mesh moving method with no
             cycle-to-cycle accumulated distortion},
  journal = {Comput. Mech.},
  volume  = {67},
  number  = {2},
  pages   = {413--434},
  year    = {2021},
  doi     = {10.1007/s00466-020-01941-y},
}

@article{DobrevTMOP2019,
  author  = {Veselin A. Dobrev and Patrick Knupp and Tzanio Kolev and Ketan Mittal and Vladimir Z. Tomov},
  title   = {The Target-Matrix Optimization Paradigm for High-Order Meshes},
  journal = {SIAM J. Sci. Comput.},
  volume  = {41},
  number  = {1},
  pages   = {B50--B68},
  year    = {2019},
  doi     = {10.1137/18M1167206},
}

@article{chiocchetti2023exactly,
  title={An exactly curl-free staggered semi-implicit finite volume scheme for a first order hyperbolic model of viscous two-phase flows with surface tension},
  author={Chiocchetti, Simone and Dumbser, Michael},
  journal={Journal of Scientific Computing},
  volume={94},
  number={1},
  pages={24},
  year={2023},
  publisher={Springer}
}

@article{dhaouadi2023structure,
  title={A structure-preserving finite volume scheme for a hyperbolic reformulation of the Navier--Stokes--Korteweg equations},
  author={Dhaouadi, Firas and Dumbser, Michael},
  journal={Mathematics},
  volume={11},
  number={4},
  pages={876},
  year={2023},
  publisher={MDPI}
}

@article{brackbill1985fluid,
  title={Fluid modeling of magnetized plasmas},
  author={Brackbill, Jeremiah},
  journal={Space Science Reviews},
  volume={42},
  number={1},
  pages={153--167},
  year={1985},
  publisher={Springer}
}

@article{brackbill1980effect,
  title={The effect of nonzero {$\nabla \cdot B$} on the numerical solution of the magnetohydrodynamic equations},
  author={Brackbill, Jeremiah U and Barnes, Daniel C},
  journal={Journal of Computational Physics},
  volume={35},
  number={3},
  pages={426--430},
  year={1980},
  publisher={Elsevier}
}

@article{dumbser2013arbitrary,
  title={On arbitrary-Lagrangian-Eulerian one-step WENO schemes for stiff hyperbolic balance laws},
  author={Dumbser, Michael and Uuriintsetseg, Ariunaa and Zanotti, Olindo},
  journal={Communications in Computational Physics},
  volume={14},
  number={2},
  pages={301--327},
  year={2013},
  publisher={Cambridge University Press}
}

@article{gaburro2026treatment,
title = {On the treatment of topology changes on 3D polyhedral moving meshes via 4D space–time hole-like elements in direct ALE ADER-DG methods},
journal = {Computer Methods in Applied Mechanics and Engineering},
volume = {463},
pages = {119434},
year = {2027},
issn = {0045-7825},
author = {Elena Gaburro and Matěj Klíma and Mauro Bonafini and Maurizio Tavelli}
}

@article{bernardelli2026structure,
  title={A structure-preserving semi-implicit finite volume scheme on vertex-staggered unstructured meshes},
  author={Bernardelli, Elena and Gaburro, Elena and Dumbser, Michael},
  journal={arXiv preprint arXiv:2604.21906},
  year={2026}
}

@article{ShashkovKnupp2002,
  author  = {Shashkov, Mikhail and Knupp, Patrick},
  title   = {Optimization-based reference-matrix rezone
             strategies for arbitrary {Lagrangian-Eulerian}
             methods on unstructured meshes},
  journal = {Sel{\c{c}}uk Journal of Applied Mathematics},
  volume  = {3},
  number  = {1},
  pages   = {81--99},
  year    = {2002}
}

@article{Loubere2010,
  author  = {Loub{\`e}re, Rapha{\"e}l and
             Maire, Pierre-Henri and
             Shashkov, Mikhail and
             Breil, J{\'e}r{\^o}me and
             Galera, St{\'e}phane},
  title   = {{ReALE}: A reconnection-based
             arbitrary-{Lagrangian--Eulerian} method},
  journal = {Journal of Computational Physics},
  volume  = {229},
  number  = {12},
  pages   = {4724--4761},
  year    = {2010},
  doi     = {10.1016/j.jcp.2010.03.011}
}

@article{Barlow2016,
  author  = {Barlow, Andrew J. and
             Maire, Pierre-Henri and
             Rider, William J. and
             Rieben, Robert N. and
             Shashkov, Mikhail J.},
  title   = {Arbitrary {Lagrangian--Eulerian} methods
             for modeling high-speed compressible
             multimaterial flows},
  journal = {Journal of Computational Physics},
  volume  = {322},
  pages   = {603--665},
  year    = {2016},
  doi     = {10.1016/j.jcp.2016.07.001}
}

@article{Sidilkover2025Spurious,
  author    = {Sidilkover, David S.},
  title     = {Spurious vorticity in Eulerian and Lagrangian methods},
  journal   = {Journal of Computational Physics},
  volume    = {520},
  pages     = {113510},
  year      = {2025},
  doi       = {10.1016/j.jcp.2024.113510}
}
\end{document}